\documentclass[11pt]{amsart}
\usepackage{amsmath,paralist,amssymb}
\usepackage{amsthm}
\usepackage{hyperref}
\usepackage{graphicx,subfigure}
\usepackage{tikz,color}
\usepackage[margin=1in]{geometry}
\usetikzlibrary{arrows,shapes,snakes,backgrounds} 
\usepackage{verbatim}

\usepackage[all]{xy}

\newcommand{\old}[1]{}
\newtheorem{theorem}{Theorem}

\newtheorem*{corollary*}{Corollary}
\newtheorem*{theorem*}{Theorem}

\newtheorem{proposition}[theorem]{Proposition}
\newtheorem{lemma}[theorem]{Lemma}
\newtheorem{corollary}[theorem]{Corollary}
\newtheorem{conjecture}[theorem]{Conjecture}
 
\theoremstyle{remark}

\theoremstyle{definition}
\newtheorem{definition}{Definition}

\newcommand{\be}{\begin{equation}}

\newcommand{\ee}{\end{equation}}

\newcommand{\bd}{\begin{definition}}
\newcommand{\ed}{\end{definition}}
\newcommand{\bt}{\begin{theorem}}
\newcommand{\et}{\end{theorem}}
\newcommand{\bl}{\begin{lemma}}
\newcommand{\el}{\end{lemma}}
\newcommand{\bp}{\begin{proposition}}
\newcommand{\ep}{\end{proposition}}
\newcommand{\bc}{\begin{corollary}}
\newcommand{\ec}{\end{corollary}}

\newcommand{\p}{{\bf p}}

\newtheorem*{lemma*}{Lemma}

\newtheorem*{theorem1*}{Theorem \ref{sp}}
\newtheorem*{theorem2*}{Theorem \ref{mm}}

\def\wt{\widetilde}

\def\i{{\rm in }}
\def\fin{{\rm fin}}

 \def\f_H{{\bf w}}
 \def\f{{\bf f}}
 
\def\R{\mathbb{R}}

\def\Z{\mathbb{Z}}

 \def\F{\mathcal{F}}
\def\O{\mathcal{O}}

\def\c{{\bf c}}

\def\alt{{\rm dep}}

 \def\x{{\bf x}}
 
 \def\b{\beta}
\def\R{\mathbb{R}}

\def\Z{\mathbb{Z}}
 \def\F{\mathcal{F}}
 
\def\C{\mathcal{C}}

\def\fin{{\rm fin}}

 \def\f_H{{\bf w}}
 \def\f{{\bf f}}
 
\def\R{\mathbb{R}}

\def\Z{\mathbb{Z}}

 \def\F{\mathcal{F}}

\def\myb{\widetilde{MACYB_n}(\mathfrak{b})}

\def\yb{\widetilde{ACYB_n}(\beta)}
\def\t{\mathcal{S}(\b)}
\def\mt{\mathcal{S}(\mathfrak{b})}
\begin{document}

\title[$h$-polynomials of reduction trees]{$h$-polynomials of reduction trees}
\author{Karola M\'esz\'aros}
\email{karola@math.cornell.edu}
\address{
Department of Mathematics, Cornell University, Ithaca, NY 14853
}
\thanks{The author was partially supported by a National Science Foundation Postdoctoral Research Fellowship  (DMS 1103933)}
\date{September 5, 2014}

\begin{abstract} 
We develop a method of proving  nonnegativity of the coefficients of certain polynomials, also called reduced forms, defined by  Kirillov in his quasi-classical Yang-Baxter algebra, its abelianization and related algebras.  It has been shown previously that the relations of the abelianization of the quasi-classical Yang-Baxter algebra, also called the subdivision algebra, encode ways of subdividing flow polytopes. In turn, these subdivisions can be represented as  reduced forms, or as reduction trees. 
We use  reduction trees in the subdivision algebra  to construct  canonical triangulations of flow polytopes which are shellable.   
We explain how a shelling  of the canonical triangulation  can be read off from the corresponding  reduction tree in the subdivision algebra.  We then introduce the notion of shellable reduction trees in the subdivision and related algebras and define $h$-polynomials of reduction trees.  In the case of the subdivision algebra, the  $h$-polynomials of  the canonical  triangulations of  flow polytopes    equal the $h$-polynomials of  the corresponding reduction trees, which motivated our definition.   We show that the reduced forms in various algebras, which can be read off from the leaves of the reduction trees,  specialize to the shifted  $h$-polynomials of the corresponding reduction trees. This yields a technique for proving nonnegativity properties of reduced forms.  As a corollary  we settle a  conjecture of A.N. Kirillov.

 
\end{abstract}

 \maketitle
\tableofcontents 

\section{Introduction}

 Nonnegativity properties abound in mathematics, and whenever one arises, the most satisfying explanation of integer nonnegativity is to demonstrate what a certain nonnegative quantity  counts. The present paper is written in this spirit and explains nonnegativity properties of polynomials using combinatorial abstractions of geometric ideas.  This is a follow up on the paper \cite{h-poly1} where the author proved that certain polynomials called shifted reduced forms in the subdivision algebra have nonnegative coefficients by showing that they equal $h$-polynomials of triangulations of flow polytopes. The methods used in \cite{h-poly1} are entirely geometric, and the purpose of this paper is the abstraction of geometric ideas related (though not identical) to those in \cite{h-poly1}. Before elaborating further, we say a few words on the reduced forms we study and their origins. 
 
 The polynomials we study are reduced forms  introduced by Kirillov in  the quasi-classical Yang-Baxter algebra and its abelianization. These algebras were defined by A.N. Kirillov \cite{Kir, k2, k2014} with Schubert calculus in mind and the former is closely related to the Fomin-Kirillov algebra \cite{fk}.  The abelianization of the  quasi-classical Yang-Baxter algebra has been considered by the present author under the name subdivision algebra, since its relations encode ways to subdivide root and flow polytopes \cite{m1, m2, m-prod}. The polynomials of interest in this paper arise as reduced forms in the above algebras; the reduced form 
 of a monomial in an algebra is obtained via substitution rules dictated by the relations of the algebra.

\textit{This paper has two components: Sections \ref{sec:shell}-\ref{sec:refining}  contain the construction and study of canonical triangulations of flow polytopes, while Sections \ref{sec:embed} and \ref{sec:conjecture7} present a method for studying reduced forms in various algebras. } These two components can be understood without reference to each other, although the ideas  in Sections \ref{sec:shell}-\ref{sec:refining}  serve as the motivation for the methods and justification for the names of the notions introduced  in Sections \ref{sec:embed} and \ref{sec:conjecture7}.
 \smallskip
 
 The essence of the subdivision algebra is that the reduced form of a monomial in it can naturally be seen as a dissection of a flow polytope corresponding to the monomial into simplices.  Any dissection obtained from a reduced form can  be graphically  represented by reduction trees,  which are simply a way of encoding a  substitution procedure dictated by the relations of an algebra.  We show that there is a \textbf{canonical triangulation} of any flow polytope that can be obtained as a dissection encoded by a reduced form (or  reduction tree), see Theorem \ref{tri}. Moreover, we show in Theorem \ref{thm-shell} that the canonical triangulations we constructed are shellable. The canonical triangulations constructed here are different from the triangulations  considered in \cite{h-poly1} and it is the geometry of the canonical triangulations that we abstract to the structural study of a reduction trees.
  
 Motivated by the considerations for  canonical triangulations,
 we establish a framework for studying reduced forms in several related algebras. This is done in Sections \ref{sec:embed} and \ref{sec:conjecture7} and we note that these sections are self-contained, and the reader interested in these applications can start there directly. We  introduce a notion alike \textbf{shellability for reduction trees} (which we call strong embeddability), inspired by the geometric notion. We  also define $h$-polynomials of reduction trees.  The \textbf{$h$-polynomials of  certain  reduction trees in the subdivision algebra   equal the $h$-polynomials of canonical  triangulations of  flow polytopes}, which was the motivation for our definition of $h$-polynomials of reduction trees.
 
  We show that the \textbf{reduced forms in various algebras specialize to the shifted  $h$-polynomi-als of the corresponding reduction trees}, see Theorems \ref{q-h} and \ref{q-ht}. This yields a technique for proving nonnegativity results  for  reduced forms.  This technique is related, though different from the one established for proving nonnegativity of shifted reduced forms in subdivision algebras in \cite{h-poly1}, since our method relies on the study of $h$-polynomials of reduction trees as opposed to $h$-polynomials of triangulations.   As a corollary to our results we settle a  conjecture of A.N. Kirillov, see Theorem \ref{7}.

The paper is organized as follows.   In Section \ref{sec:def} we  define flow polytopes. Next we explain how to subdivide flow polytopes and how we can encode the subdivisions with a reduction tree. Then we define the (multiparameter)  subdivision algebra as well as the (multiparameter) associative quasi-classical Yang-Baxter algebra of A.N. Kirillov.
 
In Section \ref{sec:shell}  we  construct canonical triangulations for flow polytopes and introduce the notion of weak  embeddability of reduction trees in order to  construct a particular shelling order for the canonical triangulations.
In  Section \ref{sec:leaves} we introduce the notion of  strong embeddability of reduction trees and indicate how to use it to give a description of the full set of leaves of the reduction tree in the special reduction order $\O$. In Section \ref{sec:refining} we study a refinement of the $h$-vector for our canonical triangulations.    

Section \ref{sec:embed}  parts from  geometry  and focuses on  the structure of reduction trees which became apparent in the previous sections. While the Sections \ref{sec:shell}-\ref{sec:refining}  are helpful for understanding the motivation for the notions in Section \ref{sec:embed}, this section is self-contained and can be read without reading the previous ones. We introduce weak and strong embeddable properties of partial reduction trees,  key notions that are seen to unify our proofs. We also define the   $h$-polynomial of a reduction tree and show that it equals the specialized (shifted) reduced form. We  generalize  our results from reduction trees to partial reduction trees. As a corollary    we prove   special cases of Conjecture 7 of A.N. Kirillov \cite{k2}  in Section \ref{sec:conjecture7} and demonstrate  via counterexamples that Conjecture 7  \cite{k2} cannot hold in its full generality. 

In Section \ref{sec:app} we prove that our canonical triangulation is indeed a triangulation. We postpone this proof to the end as it is technical, and its  ideas  are not used elsewhere in the paper.

\section{Definitions and Prelimiaries}
\label{sec:def}

For completeness,  in this section we include several key definitions used throughout the paper, following \cite{h-poly1}. For further details see \cite{h-poly1}.

\subsection{Flow polytopes and their subdivisions.} \bd \label{1} Given a loopless graph $G$  on the vertex set $[n]$,  let $\i(e)$ denote the smallest (initial) vertex of edge $e$ and $\fin(e)$ the biggest (final) vertex of edge $e$.  Let $E(G)=\{\{e_1, \ldots, e_l\}\}$   be the multiset of edges of $G$. We correspond  variables $x_{e_i}$, $i \in [l]$, to the edges of $G$, of which we think as flows.  
The {\bf flow polytope} $\F_G$ is naturally embedded into $\R^{\#E(G)}$, where  $x_{e_i}$, $i \in [l]$, are thought of as the coordinates. $\F_G$  is defined by 

 $$x_{e_i}\geq 0, \mbox{ }i \in [l],$$

  $$1=\sum_{e \in E(G), \i(e)=1}x_e= \sum_{e \in E(G),  \fin(e)=n+1}x_e,$$
  
\noindent  and for $2\leq i\leq n$
  
  $$\sum_{e \in E(G), \fin(e)=i}x_e= \sum_{e \in E(G), \i(e)=i}x_e.$$
  \ed
  \medskip
  
   Flow polytopes lend themselves to subdivisions via  {\it reductions}, as explained below. A similar property of root polytopes was established  in \cite{m2, m1}.

\bd \label{3} Given a graph $G$ on the vertex set $[n]$  containing edges $(i, j)$ and $(j, k)$, $i<j<k$, performing the {\bf reduction} on these edges of $G$ yields three graphs on the vertex set $[n]$:
  \begin{eqnarray} \label{graphs1}
E(G_1)&=&E(G)\backslash \{(j, k)\} \cup \{(i, k)\}, \nonumber \\
E(G_2)&=&E(G)\backslash \{(i, j)\} \cup \{(i, k)\},  \nonumber \\
E(G_3)&=&E(G)\backslash \{(i, j), (j, k)\} \cup \{(i, k)\}.
\end{eqnarray}

Denote by $((i,j),(j,k),L)$ the reduction that took place to get from $G$ to $G_1$,  by $((i,j),(j,k),R)$ the reduction that took place to get from $G$ to $G_2$ and $((i,j),(j,k),M)$ the reduction that took place to get from $G$ to $G_3$. $L, R, M$ correspond to ``left, right, middle."

\ed

\begin{figure}
\begin{center}
\includegraphics[scale=.7]{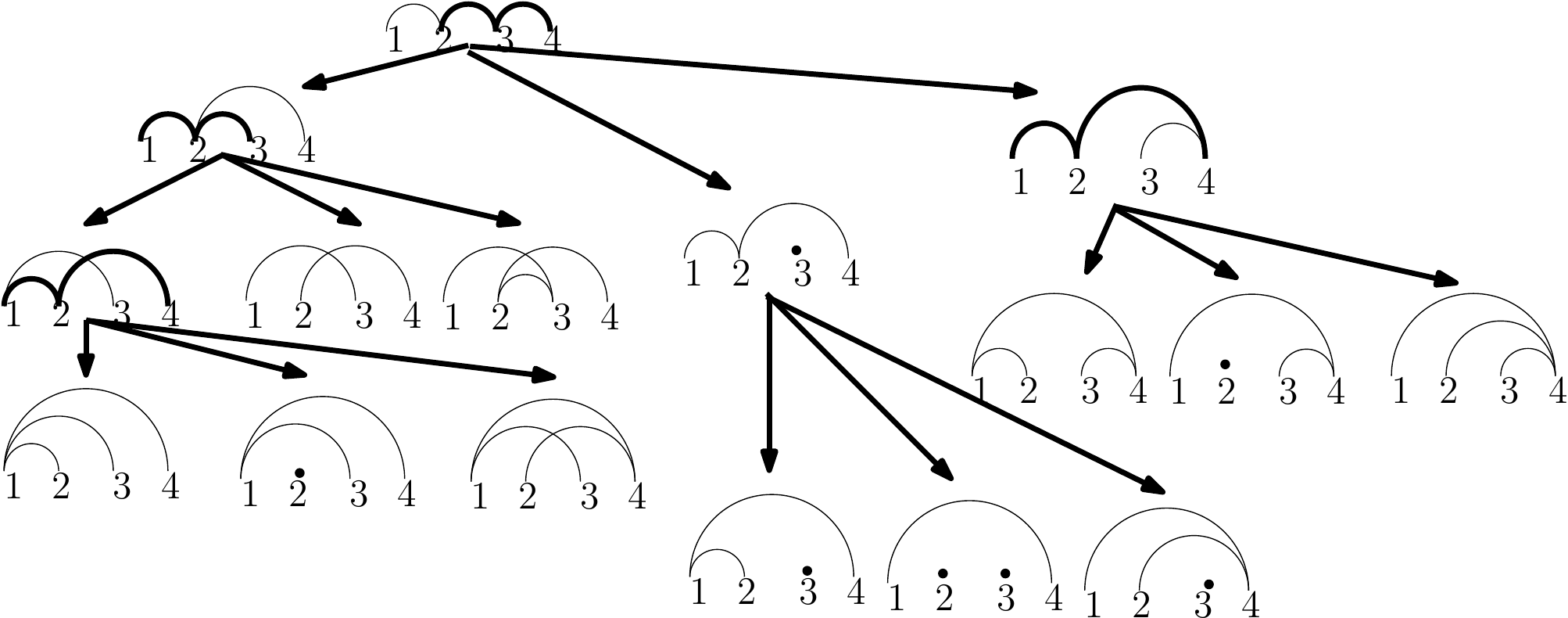}
 \caption{A reduction tree of $G=([4], \{(1,2),(2,3), (3,4)\})$. The edges on which the reductions are performed are in bold. }
 \label{fig:redtree}
 \end{center}
\end{figure}

\bd A \textbf{reduction tree} $R_G$ of a graph $G$ is a tree with nodes labeled by graphs and such that all non-leaf nodes of $R_G$ have three children. The root is labeled by $G$. If there are two edges  $(i, j), (j, k) \in E(G)$,  $i<j<k$, on which we choose to do a reduction, then the  children of the root are labeled by  $G_1, G_2$ and $G_3$ as in (\ref{graphs1}). Next, continue this way by constructing  reduction trees for $G_1, G_2$ and $G_3$. If some graph has no edges $(i, j), (j, k)$,  $i<j<k$, then it is its own reduction tree. Note that the reduction tree $R_G$ is not unique; it depends on our  choice of edges to reduce.
However, the number of leaves (referring to the graph labeling a leaf) of all reduction trees of $G$ with a given number of edges is the same, see \cite[Lemma 5]{h-poly1}. We choose a particular embedding of the reduction tree in the plane for convenience: we root it at $G$ with the tree growing downwards, and such that the left child is $G_1$, the middle child is $G_3$ and the right child is $G_2$; see Figure \ref{fig:redtree}. The leaves which have the same number of edges at the root are called \textbf{full dimensional}.
\ed

\bd \label{edges}
Let the edges of $G$ be $e_1, \ldots, e_k$, where we distinguish multiple edges. If a reduction involving edges $a=(i,j)$ and $b=(j,k)$ of $G$ is performed, then the new edge $(i,k)$ appearing in all three graphs as in \eqref{graphs1} is formally thought of as $a+b$. The other edges stay unchanged. To get to leaves $G_1$ and $G_2$  of $R_G$ we iterate this process, thereby expressing the edges of any leaf as a sum of edges of the graph being the root of the reduction tree. Two edges $c$ and $d$ in the graphs  $G_1$ and $G_2$, respectively, are the same, if they are the sum of exactly the same edges of $G$. The intersection of two graphs $G_1$ and $G_2$ in a reduction tree $R_G$ is $G_1\cap G_2=(V(G), E(G_1)\cap E(G_2))$, where if $e \in E(G_1)\cap E(G_2)$ then as explained above $e$ is the sum of the same edges of $G$ in both $G_1$ and $G_2$. 
\ed

By abuse of notation we will write $G-e$ to mean the graph $G$ with edge $e$ deleted and $G+e$ to mean the graph $G$ with edge $e$ added .

\bd   The {\bf augmented graph} ${\widetilde{G}}$ of $G=([n], E)$ is ${\widetilde{G}}=([n]\cup \{s, t\}, \widetilde{E})$, where $s$ (source) is the smallest, $t$ (target/sink) is the biggest vertex of  $[n]\cup \{s, t\}$, and $\widetilde{E}=E\cup \{(s, i), (i, t) | i \in [n]\}$.  Denote by $\mathcal{P}(\widetilde{G})$ the set of all maximal paths in $\widetilde{G}$, referred to as \textbf{routes}. It is well known that the unit flows sent along the routes in $\mathcal{P}(\widetilde{G})$ are the vertices of $\mathcal{F}(\widetilde{G})$. 
 \ed

\begin{definition} \label{leaf} Consider a node $G_1$ of the reduction tree $R_G$, where each edge of $G_1$ is considered as a sum of the edges of $G$. The image of the map $m:E(G_1)\rightarrow \mathcal{P}(\widetilde{G})$ which takes an edge $(v_1, v_2)=e=e_{i_1}+\cdots+e_{i_l}$, $e \in G_1$, $e_{i_j} \in E(G)$, $j \in [l]$, to the route $(s, v_1), e_{i_1},\ldots,e_{i_l}, (v_2, t)$ gives the vertices of $\F_{\widetilde{G_1}}$ (by taking the unit flows  on these routes).  In case $G_1$ is not a node of the reduction tree $R_G$, but it is  an intersection of  nodes of $R_G$, so that each edge of $G_1$ can still be considered as a sum of the edges of $G$, we still define $\F_{\widetilde{G_1}}$ as above. This definition of  $\F_{\widetilde{G_1}}$ is of course with respect to $G$, and this is understood from the context. 
\end{definition}

Using the above definitions the proof of  the following lemma is an easy exercise.

\bl{ \cite[Proposition 1]{m-prod},\cite[Proposition 4.1]{mm}, \cite{p, S}}
 \label{red} Given a graph $G$ on the vertex set $[n]$ and   $(i, j), (j, k) \in E(G)$,   for some $i<j<k$, with $G_1, G_2, G_3$ as in \eqref{graphs1}  and $ \F_{{\widetilde{G_i}}}$, $i \in [3]$, as in Definition \ref{leaf} we have  $$\F_{{\widetilde{G}}}=\F_{{\widetilde{G_1}}} \bigcup \F_{{\widetilde{G_2}}},  \F_{{\widetilde{G_1}}} \bigcap \F_{{\widetilde{G_2}}}=\F_{{\widetilde{G_3}}} \text{ and } \F_{{\widetilde{G_1}}}^\circ \bigcap \F_{{\widetilde{G_2}}}^\circ=\emptyset,$$

 \noindent where $\F_{{\widetilde{G}}}$, $\F_{{\widetilde{G_1}}}$, $\F_{{\widetilde{G_2}}}$ are of the same dimension $d-1$, $\F_{{\widetilde{G_3}}}$ is $d-2$ dimensional,  and  $\mathcal{P}^\circ$ denotes the interior of $\mathcal{P}$.
\el

\subsection{Algebras related to flow polytopes.} Note that the reduction of graphs given in \eqref{graphs1} can be encoded as the following relation:

\be x_{ij}x_{jk}=x_{ik}x_{ij}+x_{jk}x_{ik}+\b x_{ik},  \mbox{ for }1\leq i<j<k\leq n. \ee

If we wanted to preserve more information on the actual reduction, we could instead consider the following relation: 

\be x_{ij}x_{jk}=x_{ik}x_{ij}+x_{jk}x_{ik}+\b_i x_{ik},  \mbox{ for }1\leq i<j<k\leq n. \ee

These relations give rise to what we call {\it subdivision algebras}. 

\bd The multiparameter associative \textbf{subdivision algebra} of weight $\mathfrak{b}=(\b_1, \ldots, \b_{n-1})$, denoted by $\mt$, is an associative algebra, over the ring of polynomials $\Z[\b_1, \ldots, \b_{n-1}]$, generated by the set of elements $\{x_{ij} : 1 \leq i<j \leq n\}$, subject to the relations:

(a) $x_{ij}x_{kl}=x_{kl}x_{ij},$ if  $i<j$, $k<l$,

(b) $x_{ij}x_{jk}=x_{ik}x_{ij}+x_{jk}x_{ik}+\b_i x_{ik},$ if $1\leq i<j<k\leq n$.

Letting $\b_i=\b$, $ i \in [n-1]$, the algebra $\mt$ specializes to the subdivision algebra of weight $\b$ denoted by $\t$.
\ed

The algebra $\t$ has been studied  in \cite{m2}  and  \cite{h-poly1}.

\bd \label{rf} Given a monomial $M$  in $\t$ or $\mt$, its reduced form is defined as follows.  Starting with $p_0=M$, produce a sequence of polynomials $p_0, p_1, \ldots, p_m$ in the following fashion.  To obtain $p_{r+1}$ from $p_r$,  choose a term of   $p_r$ which  is divisible by $x_{ij}x_{jk}$, for some $i,j,k$, and replace the factor  $x_{ij}x_{jk}$ in this term   with   $x_{ik}x_{ij}+x_{jk}x_{ik}+\b x_{ik}$ or $x_{ij}x_{jk}=x_{ik}x_{ij}+x_{jk}x_{ik}+\b_i x_{ik}$, depending on which algebra we are in. Note that   $p_{r+1}$   has two more terms than $p_r$. Continue  this process until   a  polynomial  $p_m$  is obtained, in which  no term is divisible by  $x_{ij}x_{jk}$, for any $i,j,k$.  Such a polynomial $p_m$  is a {\bf reduced form} of $M$. Note that we allow the use of the commutation relations of each algebra in this process.
\ed

Given a monomial $M$ in $\t$ or $\mt$  we can encode it by a graph $G_M$, simply by letting the edges of $G$ be the given by the indices of the variables in $M$. Denote a reduced form of $M$  in $\t$ by $Q_{G_M}^{\t}(\x; \b)$ and the reduced form  of $M$  in $\mt$ by $Q_{G_M}^{\mt}(\x; \mathfrak{b})$. If in the reduced forms we set $\x=(1, \ldots, 1)$, then in the notation we omit $\x$:  
$Q_{G_M}^{\t}(\b)$ or $Q_{G_M}^{\mt}(\mathfrak{b})$.

It is easy to see that by definition, the reduced form of a monomial in the subdivision algebras can be read off from the reduction tree of the corresponding graph obtained by simply taking its edge set to the the double indices of the variables of the monomial. 

Note  that the reduced form of a monomial in $\mt$ or $\t$ is not necessarily unique, which could be a desirable property. The noncommutative counterpart of $\t$, denoted by $\yb$ and defined by Kirillov \cite{k2, k2014},  is much like $\t$, but with reduced forms  unique \cite{m2}. While the same is not true of the similar noncommutative generalization, denoted $\myb$,  of $\mt$, this algebra also has beautiful combinatorics. It was A.N. Kirillov \cite{k2, k2014} who introduced these algebras and shed the first light on their rich combinatorial structure. 

\bd \cite[Definitions 3.1 and 3.2]{k2} The multiparameter associative \textbf{quasi-classical Yang-Baxter algebra} of weight $\mathfrak{b}=(\b_1, \ldots, \b_{n-1})$, denoted by $\myb$, is an associative algebra, over the ring of polynomials $\Z[\b_1, \ldots, \b_{n-1}]$, generated by the set of elements $\{x_{ij} : 1 \leq i<j \leq n\}$, subject to the relations:

(a) $x_{ij}x_{kl}=x_{kl}x_{ij},$ if  $\{i, j\}\cap \{k, l\}=\emptyset$,

(b) $x_{ij}x_{jk}=x_{ik}x_{ij}+x_{jk}x_{ik}+\b_i x_{ik},$ if $1\leq i<j<k\leq n$.

Letting $\b_i=\b$, $ i \in [n-1]$, the algebra $\myb$ specializes to the  associative quasi-classical Yang-Baxter algebra of weight $\b$, denoted by $\yb$. 

\ed

The definition of reduced forms in $\yb$ and $\myb$ is the similar to Definition \ref{rf}; the only difference is that now the order of variables matters and so we take consecutive variables $x_{ij}$ and $x_{jk}$ and replace them by $x_{ik}x_{ij}+x_{jk}x_{ik}+\b_i x_{ik}.$ E.g., while in $\mt$ the monomial $x_{jk}x_{ij}$ could be reduced, it is itself a reduced form in $\myb$.  One can also define the notion of a reduction tree for these algebras, which, in the terminology of Section \ref{sec:embed} can be seen as partial reduction trees with respect to the subdivision algebra.

\section{Weak embeddability and shelling triangulations}
\label{sec:shell}

In \cite{h-poly1} we studied several regular triangulations of $\F_{\widetilde{G}}$ relying on the work of Danilov, Karzanov and Koshevoy \cite{kosh}. They have posed the question of whether the triangulations they construct in \cite{kosh} are all of the regular triangulations of flow polytopes. The aims of this section are twofold. First, we construct a triangulation of $\F_{\widetilde{G}}$ which is not one of those constructed in \cite{kosh}. We prove that our triangulation is shellable. We leave the question of whether it is regular open for further investigation. Second,  we introduce  the notion of \textit{weak embeddabiliy of reduction trees} which   can be extended even when we are not using the   geometry of triangulations.   

\begin{figure}
\begin{center}
\includegraphics[scale=.7]{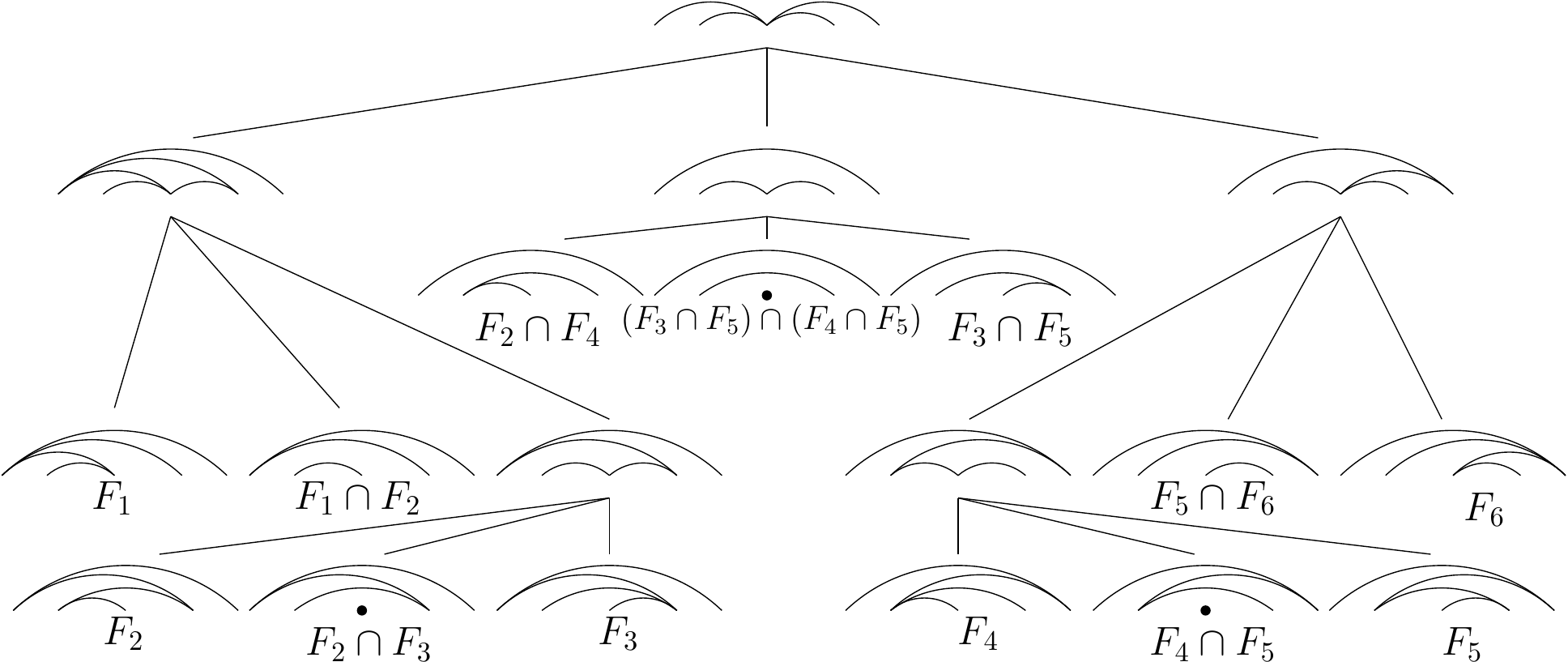}
 \caption{The reduction tree of $G=([5], \{(1,3),(2,3), (3,4), (3,5)\})$ with reductions executed in order $\O$.  The labels $F_i$, $i \in [6]$,  are explained in Theorems \ref{thm-shell} and \ref{thm-shell2}.}
 \label{fig:o}
 \end{center}
\end{figure}

The triangulation we consider in this section is obtained from the reduction tree $R_G^{\O}$, which is a reduction tree where the reductions are executed in a certain order ${\O}$. The reduction order $\O$ is defined as follows.  Given an arbitrary graph $G$ on the vertex set $[n]$, do the reductions in $G$ proceeding from the smallest vertex towards the greatest in order. Look for the smallest vertex $v$ which is nonalternating, that is that has both an edge $(a, v)$ and an edge $(v, b)$ incident to it, with $a<v<b$. Look at the two topmost edges at $v$, that is edges  $(a, v)$  and  $(v, b)$ such that $a<v<b$ and there are no edges $(a', v)$ with $a'<a$ and $(v, b')$ with $b<b'$. Do the reduction on the two topmost edges at $v$. Continue in this fashion on each leaf of the partial reduction tree ultimately arriving to the reduction tree $R_G^{\O}$ with all leaves alternating graphs, that is all of their vertices are alternating. For a  reduction tree $R_G^{\O}$ see Figure \ref{fig:o}.
 
Since the proof is technical, and not of central importance to the rest of the  paper, we state Theorem \ref{tri} here and refer the reader to Section \ref{sec:app} for a proof. 

\begin{theorem} \label{tri} The simplices corresponding to the  full dimensional leaves of  $R_{G}^{\O}$ induce a triangulation; that is, the intersection of any two of them is a face of both. Moreover, the simplices corresponding to all leaves of  $R_{G}^{\O}$ are part of this triangulation. \end{theorem}

We note that the  set of $R_{G}^{\O}$-triangulation we obtain as described in Theorem \ref{tri} are not a subset  of the  triangulations constructed in \cite{kosh}; for example, considering the graph $G=([6], \{(1,3),(3,4),$ $(4,5),(2,4),(4,6)\})$, regardless of the framing of $\widetilde{G}$, we can obtain routes which are noncoherent at vertex $4$ to be vertices of a top dimensional simplex in  the  $R_{G}^{\O}$-triangulation.

Instrumental in this section is the order of the leaves of $R_G^{\O}$: let $F_1, \ldots, F_l$  be the full dimensional leaves of  $R_G^{\O}$  in depth-first search order as shown in Figure \ref{fig:o}. Remember that we have an embedding of $R_G^{\O}$ in the plane where $G$ is the root and the graphs $G_1, G_2, G_3$ as in \eqref{graphs1} are the left, right, middle child, respectively. Also, by Theorem \ref{tri} the simplices   $\F_{\widetilde{F_1}}, \ldots, \F_{\widetilde{F_l}}$ are the top dimensional simplices in a triangulation of $\F_{\widetilde{G}}$; we refer to this triangulation as the \textbf{canonical triangulation} of $\F_{\widetilde{G}}$. Theorem \ref{thm-shell}  is the main result of  this section, and the weak embeddable property in Definition \ref{weak} introduced to prove this theorem is the essential ingredient we carry forward to less geometric settings.

\bt \label{thm-shell}  $\F_{\widetilde{F_1}}, \ldots, \F_{\widetilde{F_l}}$ is a shelling order of the canonical triangulation of $\F_{\widetilde{G}}$. \et 


\begin{definition} \label{weak} A reduction tree $R_G$ is said to have the {\bf (right) weak embeddable property} if one of the following is true for every node $H$ of $R_G$:  

1. $H$ is a leaf

2. the middle child of $H$ is $H_3$ and the right child of $H$ is $H_2$, satisfying that there is a map $b_H$ from the full dimensional leaves of the subtree $R_{H_3}$ of $R_G$ (leaves with $|E(H_3)|$ number of edges) into  the full dimensional leaves of the subtree $R_{H_2}$ (leaves with $|E(H_2)|$ number of edges) of $R_G$ such that if $b_H(L)=L'$, then $E(L')=E(L) \cup \{e\}$ with $e \not \in E(L)$. Moreover, if $L'$ is a full dimensional leaf of $R_{H_2}$ such that there is a leaf $L$ of  $R_{H_3}$ with the property that $E(L')=E(L) \cup \{e\}$ with $e \not \in E(L)$, then $L$ is in the image of $b_H$ and $b_H$ is a bijection from the full dimensional leaves of  $R_{H_3}$  onto its image. Moreover, there is a unique $L$ in $R_{H_3}$ such that  $E(L')=E(L) \cup \{e\}$ with $e \not \in E(L)$.
\end{definition}


\begin{definition} \label{leftweak} A reduction tree $R_G$ is said to have the {\bf left weak embeddable property} if it satisfies the conditions of Definition \ref{weak} when we replace $H_2$ by $H_1$ in the statement. 
\ed

\bd  \label{twosidedweak} A reduction tree $R_G$ is said to have the {\bf twosided weak embeddable property} if it has both the right and the left weak embeddable property. 
\ed

\begin{lemma} \label{embed} The reduction tree $R_G^{\O}$ has the twosided weak embeddable property.
\end{lemma}

Before proving Lemma \ref{embed}, we define the map $b_H$ on its non-leaves which we will show satisfies 2. in Definition \ref{weak}.

\bd When performing a reduction $(a, b, X)$, $X \in \{M, L, R\}$ (notation as  in Definition \ref{3}), we say that the edge $a$ is {\bf dropped} when $X=M, R$, and $a$ is {\bf kept} if $X=L$. Similarly, edge $b$ is dropped when $X=M, L$, and $b$ is kept if $X=R$. We also say that an edge $e$  is {\bf derived from the edge $b$} if it resulted as a sequence of reductions involving $b$, or sums of edges of $b$, and $e$ itself is a sum of edges with $b$.  We also signal this by saying that $e$ is a $b^*$-edge.
\ed

\bd \label{bh} Consider a non-leaf node $H$ of  $R_G^{\O}$ and let the  reduction  performed  at $H$ be  $(a,b)$ yielding middle child $H_3$ and right child $H_2$. 
Define the map $b_H$ from the full dimensional leaves of $R_{H_3}^{\O}$, denoted $FL_{H_3}^{\O}$, into the full dimensional leaves of $R_{H_2}^{\O}$,  denoted $FL_{H_2}^{\O}$, as follows. Let $(a,b, M), (c_1, d_1, X_1), (c_2, d_2, X_2), \ldots, (c_k, d_k, X_k)$ for edges $c_i, d_i$, and $X_i \in \{L, R\}$, $i \in [k]$, be the sequence of reductions leading to $L \in FL_{H_3}^{\O}$ from $H$. Recall that since the order $\O$ of reductions is specified, we are only wondering at each step whether to go $L$ or $R$. Then $b_H(L)$ is defined as the element of   $FL_{H_2}^{\O}$ obtained from $H$ by the sequence of reductions $(a, b, R), (a_1, b_1, Y_1), (a_2, b_2, Y_2), \ldots, (a_l, b_l, Y_l)$ (here again the pair of edge $a_i, b_i$ is determined by the sequence $(a, b, R), (a_1, b_1, Y_1), (a_2, b_2, Y_2), \ldots, (a_{i-1}, b_{i-1}, Y_{i-1})$), where if $(a_i, b_i)=(c_j, d_j)$, then $Y_j=X_j$ and if $(a_i, b_i)\neq (c_j, d_j)$, then it follows that $a_i$ or $b_i$ is an edge derived from $b$, but not derived from $a+b$. In this case we choose $Y_i$ so that the edge derived from $b$ is dropped in the reduction. \ed

\noindent{\it Proof of Lemma \ref{embed}.} Consider a node $H$ of  $R_G^{\O}$. If $H$ is a leaf, there is nothing to check. If $H$ is not a leaf, we claim that the map $b_H$ defined in  Definition \ref{bh} satisfies property 2. of Definition \ref{weak}.  In this proof when we refer to a $b^*$-edge, we mean a $b^*$-edge not derived from $a+b$. 

Let $S$ be the set of full dimensional leaves of $R_{H_2}^{\O}$ obtained by a sequence of  reductions that either do not involve a $b^*$-edge, or if the reduction involves a $b^*$-edge, then in the reduction we go towards the outcome where this $b^*$-edge is dropped. Clearly, the image of $b_H$ is in $S$. We now show that the inverse of $b_H$ is defined on $S$. Indeed, if the leaf $L'$ in $S$ was obtained by a series of reductions on the pairs of edges $(a, b, R), (a_1, b_1, X_1), (a_2, b_2, X_2), \ldots, (a_l, b_l, X_l)$ of which  $(a_{i_1}, b_{i_1}), (a_{i_2}, b_{i_2}), \ldots, (a_{i_m}, b_{i_m})$ are the ones involving $b^*$-edges, then the sequence of reductions  $(a, b, M), (a_1, b_1, X_1), (a_2, b_2, X_2), \ldots, (a_l, b_l, X_l)$ with the reductions   $(a_{i_1}, b_{i_1}, X_{i_1}), (a_{i_2}, b_{i_2}, X_{i_2}), \ldots,$  $(a_{i_m}, b_{i_m}, X_{i_m})$ deleted, is a valid sequence of reductions in $R_{H_3}^{\O}$ leading to a leaf $L$, and the map $b_H$ takes $L$ to $L'$. Note also that if $L' \not \in S$ is a full dimensional leaf of $R_{H_2}^{\O}$ not in $S$, then it contains more than one $b^*$-edge, in which case it cannot be that $E(L')=E(L) \cup \{e\}$ with $e \not \in E(L)$, for some leaf $L \in FL_{H_3}^{\O}$. Moreover, if $L' \in S$, then it has a unique $b^*$-edge and thus it is the image of a unique $L$ in $R_{H_3}$. 

To prove that $R_G^{\O}$ has the left embedabble property, one can define an analogous map $b_H^L$, where the role of $b$ is played by the edge $a$. 
\qed

\medskip

 Next we define the depth of a reduction tree. 
 
\bd \label{dep} 
Let $\alt(G)$, the \textbf{depth of $R_G^{\O}$},  be  the maximum length of a path in $R_G^{\O}$ from $G$ to a leaf. 
\ed

\noindent \textit{Proof  of Theorem \ref{thm-shell}.}     We prove by induction on $\alt(G)$,  that $\F_{\widetilde{F_1}}, \ldots, \F_{\widetilde{F_l}}$ is a shelling order.  Let $G_1, G_2,$ and  $G_3$ be the left, right, and middle child of $G$, respectively, in $R_G^{\O}$ after having performed a reduction on edges $a$ and $b$. Let $F_1, \ldots, F_k$ be the full dimensional leaves of $R_{G_1}^{\O}$,  $F_{k+1}, \ldots, F_l$ be the full dimensional leaves of $R_{G_2}^{\O}$, and  $Q_1, \ldots, Q_z$ be the full dimensional leaves of $R_{G_3}^{\O}$ in depth first search order.  By induction hypothesis  $F_1, \ldots, F_k$ and  $F_{k+1}, \ldots, F_l$ and $Q_1, \ldots, Q_z$  are shelling orders of the canonical triangulations of $\F_{\widetilde{G_1}}, \F_{\widetilde{G_2}}$ and $\F_{\widetilde{G_3}}$, respectively, obtained via the order $\O$.  There are several things we  prove in order to prove that $\F_{\widetilde{F_1}}, \ldots, \F_{\widetilde{F_l}}$ is a shelling order. Let $S$ be the image of $b_G$  in $R_{G_2}^{\O}$ as constructed in the proof of Lemma \ref{embed} and let  $\bar{S}$ be the full dimensional leaves of $R_{G_2}^{\O}$ not in $S$. Then, to prove Theorem \ref{thm-shell} it suffices to prove Claims 1 and 2:  

{\bf Claim 1.}  for $F_i \in \bar{S}$,  \be \F_{\widetilde{F_i}} \cap (\F_{\widetilde{F_1}} \cup \cdots \cup \F_{\widetilde{F_{i-1}}})=  \F_{\widetilde{F_i}} \cap (\F_{\widetilde{F_{k+1}}} \cup \cdots \cup \F_{\widetilde{F_{i-1}}})\ee

{\bf Claim 2.} for $F_i \in {S}$,  \be \F_{\widetilde{F_i}} \cap (\F_{\widetilde{F_1}} \cup \cdots \cup \F_{\widetilde{F_{i-1}}})= \big( \F_{\widetilde{F_i}} \cap (\F_{\widetilde{F_{k+1}}} \cup \cdots \cup \F_{\widetilde{F_{i-1}}}) \big) \cup \F_{\widetilde{b_G^{-1}(F_i)}}\ee
 

By Lemma \ref{red} we have \be \F_{\widetilde{G_1}} \cap \F_{\widetilde{G_2}}=\F_{\widetilde{G_3}}, \ee which can also be written as

 \be \label{cond} ( \F_{\widetilde{F_1}} \cup \cdots \cup  \F_{\widetilde{F_k}})\cap (  \F_{\widetilde{F_{k+1}}} \cup \cdots \cup  \F_{\widetilde{F_l}})=  \F_{\widetilde{Q_1}} \cup \cdots \cup  \F_{\widetilde{Q_z}}. \ee

 Let $b_G(Q_j)=F_{i_j}$, $j \in [z]$, $k+1\leq i_j \leq l$; that is, $F_{i_j}=Q_j+e$.  Using Theorem \ref{coord} 
for the graphs $Q_j$ and $F_{i_j}$ together with Corollary \ref{dimension} 
  we can conclude that $\F_{\widetilde{Q_j}}$ is a facet of $\F_{\widetilde{F_{i_j}}}$.  Using in addition the properties of $b_G$ as given in Definition \ref{weak}  we can also conclude that   $\F_{\widetilde{Q_j}}$ is not a facet of any other $\F_{\widetilde{F_{i}}}$  for $k+1\leq i \leq l$, $i \neq i_j$. Moreover, since   $\F_{\widetilde{F_1}}, \ldots, \F_{\widetilde{F_l}}$  are the top dimensional simplices in the canonical triangulation of $ \F_{\widetilde{G}}$ for which \eqref{cond} holds, we have that   \be \label{}  \F_{\widetilde{Q_{j}}}\subset (\F_{\widetilde{F_1}} \cup \cdots \cup \F_{\widetilde{F_k}})\cap  \F_{\widetilde{F_{i_j}}}. \ee 

By the above together with Theorem \ref{tri},
 we have that if $\F_{\widetilde{F_i}}$, $k+1\leq i \leq l$, is not in the image of $b_G$ then $\F_{\widetilde{F_i}}$ does not attach on a facet to $(\F_{\widetilde{F_1}} \cup \cdots \cup \F_{\widetilde{F_k}})$, and if $\F_{\widetilde{F_i}}$, $k+1\leq i \leq l$, is  in the image of $b_G$ then $\F_{\widetilde{F_i}}$ attaches on exactly one facet $\F_{\widetilde{b_G^{-1}(F_i)}}$  to $(\F_{\widetilde{F_1}} \cup \cdots \cup \F_{\widetilde{F_k}})$. Thus,   we have identified all the facets on which the  $\F_{\widetilde{F_i}}$'s, $k+1\leq i \leq l$, attach to previous simplices in the canonical triangulation and they agree with the facets specified in Claims 1 and 2 above. In order to finish the proof of Claims 1 and 2, and thus that we have a shelling, it remains to prove that the   $\F_{\widetilde{F_i}}$'s, $k+1\leq i \leq l$, only attach on facets and not on lower dimensional faces to 
$(\F_{\widetilde{F_1}} \cup \cdots \cup \F_{\widetilde{F_{i-1}}})$. This is what we do next.

  In light of Theorem \ref{coord}, Corollary \ref{dimension} and Theorem \ref{subsets} 
  what needs to be checked is as follows: if $H_a=F_a\cap F_i$, for $a \in [k]$ and some fixed $i \in [k+1, l]$, with $|E(H_a)|<|E(F_i)|-1$,  then $H_a \subset F_i \cap b_G^{-1}(F_i)$, if this is well defined, or $H_a \subset   (F_r\cap F_i)$, for some $r \in [k+1, i-1]$. Note that if $b_G^{-1}(F_i)$ is well defined, then $F_i$ has exactly one   $b^*$-edge $e$, which is not derived from $a+b$, where the reduction at $G$ is performed on the edges $a$ and $b$,   and this $b^*$-edge $e$ cannot appear in any leaf of $R_{G_3}^{\O}$, thus $H_a \subset  b_G^{-1}(F_i)=F_i-e$. Also note that if $\F_{\widetilde{F_i}}$ attaches on at least two facets to $(\F_{\widetilde{F_{k+1}}} \cup \cdots \cup \F_{\widetilde{F_{i-1}}})$, say as an intersection with $\F_{\widetilde{F_c}}$ and $\F_{\widetilde{F_d}}$, $k+1\leq c\neq d\leq i-1$, then $H_a$ is subset of at least one of $F_c$ or $F_d$. Thus, it remains to deal with the case where $b_G^{-1}(F_i)$ is not well defined (and so $F_i$ has at least two $b^*$-edges not derived from $a+b$) and $F_i$ attaches on exactly one facet to $(\F_{\widetilde{F_{k+1}}} \cup \cdots \cup \F_{\widetilde{F_{i-1}}})$. Since  $b_G^{-1}(F_{k+1})=Q_1$, then $i>k+1$. Since $\F_{\widetilde{F_{k+1}}}, \ldots, \F_{\widetilde{F_{l}}}$ is a shelling order by induction, it follows that $\F_{\widetilde{F_i}}$ attaches on exactly one facet to $(\F_{\widetilde{F_{k+1}}}\cup \cdots \cup \F_{\widetilde{F_{i-1}}})$ and on no lower dimensional face.

 Let the facet on which $\F_{\widetilde{F_i}}$ attaches to $(\F_{\widetilde{F_{k+1}}}\cup \cdots \cup \F_{\widetilde{F_{i-1}}})$  be the intersection of $\F_{\widetilde{F_i}}$ and $\F_{\widetilde{F_j}}$, $k+1\leq j<i$,  with $F_i \cap F_j=F_i-f(d)$ for some edge $f(d)$ of $F_i$, which we describe shortly. Consider the path from $G$ to $F_i$ in  $R_G^{\O}$. Let the node of $R_G^{\O}$ where the last right edge is taken on the path from $G$ to $F_i$ be $H$. Let the reduction done at $H$ be on the edges $c=(i,j)$ and $d=(j,k)$ for $i<j<k$. Let $H_1, H_2, H_3$ be the left, right and middle children of $H$. Then, $F_i=b_H(Q)$ for a graph $Q \in R_{H_3}^{\O}$, where $E(F_i)=E(Q) \cup \{f(d)\}$, where $f(d)$ is a $d^*$-edge not derived from $c+d$. 
  Since $\F_{\widetilde{F_i}}$ attaches on exactly one facet to $(\F_{\widetilde{F_{k+1}}}\cup \cdots \cup \F_{\widetilde{F_{i-1}}})$ and on no lower dimensional face, it follows that the graphs $F_{k+1}, \ldots, F_{i-1}$ do not contain the edge $f(d)$. 

Let the sequence of reductions leading from $G$ to $F_i$ be $s=(s_0, s_1, \ldots, s_z)$, where $s_0=(a, b, R)$ and $s_i=(c_i, d_i, X_i)$ for some pair of edges $c_i, d_i$, $i \in [z]$, and $X_i \in \{L, R\}$. We now establish what it means for edge $f(d)$ to be an edge of a leaf of $R_G^{\O}$, in terms of the reductions leading to it. 

Let $f(d)=e_1+\ldots+e_l$, where $e_i$, $i \in [l]$, are edges of $G$. Let the reductions on $e_1$ and $e_2$, on $e_1+e_2$ and $e_3, \ldots,$ on $e_1+\ldots+e_{l-1}$ and $e_l$ be subsequence of $s$ (with appropriate ordering among the edges of the pairs and L and R added). Denote this subsequence of $s$ by $s^e$, where the superscript $e$  signifies that these reductions are essential in creating $f(d)$. Note that in order to be able to do these reductions, it is also key that the edges we want to do the reduction on are present, that is, if $e_1+\ldots+e_i$, $i \in [l]$, is an edge which is part of the reduction with an edge other than $e_{i+1}$, we must keep it in that reduction. 
Since the order of reductions is prescribed by $\O$ it follows that  the edge $f(d)$ is an edge of a leaf of $R_G^{\O}$ if and only if $s^e$ is a subsequence of the sequence of reductions leading to that leaf.  

There cannot be a graph $H'$ in $R_{G_2}^{\O}$ preceeding $F_i$ (meaning that the path from $G_2$ to $H'$ is to the left of the path from $G_2$ to $F_i$) such that $s^e$ is a subsequence of   the sequence of reductions leading from $G$ to  $H'$, since then a descendent of $H'$ would contain $f(d)$. 
We need to prove using this and that $b_G^{-1}(F_i)$ is not defined that there is also no graph $H'$ in $R_{G_1}^{\O}$  such that $s^e$ is a subsequence of   the sequence of reductions leading from $G$ to  $H'$. 

Since $b_G^{-1}(F_i)$ is not defined, it follows that $s$ contains a reduction involving the edge $b$ (other than $(a, b, R)$) and moreover, it also contains a reduction involving a $b^*$-edge $e$ not derived from $a+b$ where that edge $e$ is kept after the reduction is performed, thereby creating at least two $b^*$-edges not derived from $a+b$. Obviously, if any of the reductions involving $b^*$-edges   are among $s^e$, then $f(d)$ cannot appear in $R_{G_1}^{\O}$.
We argue that if none of the reductions involving $b^*$-edges  are among $s^e$, then there is a graph $H'$ in $R_{G_2}^{\O}$ preceeding $F_i$ such that $s^e$ is a subsequence of   the sequence of reductions leading from $G$ to  $H'$, which would contradict our assumption that $\F_{\widetilde{F_i}}$ attaches on exactly one facet to $(\F_{\widetilde{F_{k+1}}}\cup \cdots \cup \F_{\widetilde{F_{i-1}}})$ and on no lower dimensional face.

We  now elaborate why under the above circumstances if none of the reductions involving $b^*$-edges not derived from $a+b$ are among $s^e$, then there is a graph $H'$ in $R_{G_2}^{\O}$ preceeding $F_i$ such that $s^e$ is a subsequence of   the sequence of reductions leading from $G$ to  $H'$. Since $s$ contains a reduction involving the edge $b$ (other than $(a, b, R)$) and  it also contains a reduction involving a $b^*$-edge not derived from $a+b$ and where that edge is kept, it follows that for some edge $x$ there is a reduction $(x, b, R)$ in $s$  and $x\neq e_1+\cdots+e_i$, for any $i \in [l]$, or there is a reduction $(g(a+b), x, R)$  in $s$ where $g(a+b)$ is an $(a+b)^*$-edge and  $x\neq e_1+\cdots+e_i$, for any $i \in [l]$, and this reduction is followed by a reduction $(f(b), x, L)$, where $f(b)$ is a $b^*$-edge not derived from $a+b$.  However, if none of the reductions involving $b^*$-edges not derived from $a+b$ are among $s^e$,  then there is a graph $H'$ in $R_{G_2}^{\O}$ in the subtree to which we get if we do  $(x, b, L)$ instead of $(x, b, R)$ or $(g(a+b), x, L)$ instead of $(g(a+b), x, R)$,  such that $s^e$ is a subsequence of   the sequence of reductions leading from $G$ to  $H'$. 
\qed

 \section{Strong embeddability and a description of the leaves of $R_G^{\O}$}
\label{sec:leaves}

In this section we introduce \textit{strong embeddability of reduction trees} and  use it to give a description of all the leaves in  $R_G^{\O}$. As we will see in in Sections \ref{sec:embed} and \ref{sec:conjecture7}  strong embeddability generalizes to other settings. 

\bt \label{thm-shell2}   Let ${F_1}, \ldots, {F_l}$ be the full dimensional leaves of $R_G^{\O}$ in depth-first search order. Let $$P_i:=\{\{Q_1^i, \ldots, Q_{f(i)}^i\}\}=\{\{F_i \cap F_j \mid 1\leq j<i, |E(F_i \cap F_j)|=|E(F_i)|-1\}\}.$$  Then 
\be \label{formal} \sum_{i=1}^l \prod_{j=1}^{f(i)} (F_i+Q_j^i)\ee

\noindent is the formal sum of the  set of the leaves of $R_G^{\O}$, where the product of graphs is their intersection, and if $f(i)=0$ then we define $\prod_{j=1}^{f(i)} (F_i+Q_j^i)=F_i$.

\et

Note that with the notation of Theorem \ref{thm-shell2} and with Definition \ref{leaf} in mind we have that  

\be \F_{\widetilde{F_i}} \cap (\F_{\widetilde{F_1}}  \cup \cdots \cup \F_{\widetilde{F_{i-1}}})=\F_{\widetilde{Q^i_1}} \cup  \cdots \cup \F_{\widetilde{Q^i_{f(i)}}},\nonumber \ee where $\F_{\widetilde{Q^i_j}}$, $j \in [f(i)]$, is a facet of $\F_{\widetilde{F_i}}$. Indeed, this follows directly from Theorems \ref{coord} 
and Theorem \ref{thm-shell}.

\bd \label{facet} Given a reduction tree  $R_G$ and a full dimensional leaf $L$ of it, we say that a leaf $H$ of  $R_G$  is a {\bf preceeding  facet} of $L$ if 

1. $H$ is before $L$ in the depth first search order of  the leaves of   $R_G$ 

2. $E(H)\subset E(L)$ and $|E(H)|=|E(L)|-1$

3. the unique path in  $R_G$ from $L$ to $H$ consists of several up steps followed by several down steps, so that the first of the down steps is a Middle reduction. 

\ed

\bl  \label{q_i} Let ${F_1}, \ldots,{F_l}$ be full dimensional leaves of $R_G^{\O}$ in depth-first search order  and let $P_i$ be as in Theorem \ref{thm-shell2}. Then the (multi)set of preceeding facets of $F_i$ is equal to $P_i$.
\el

Before proving Lemma \ref{q_i} we provide an auxiliary lemma that will come in handy in the proof.

\bl \label{differ} Let $F_i$ and $F_j$, with $F_j$ preceeding $F_i$ in depth-first search order, be two full dimensional leaves of $R_G^{\O}$ differing in only one edge; that is, $F_i\cap F_j= F_i-e$ for some edge $e$ of $F_i$. Let the paths from $G$ to $F_j$ and from $G$ to $F_i$ split at graph $H$ via the reduction $(z, v)$, where we take $(z, v, L)$ towards $F_j$ and $(z, v, R)$ towards $F_i$. Then:

1. on the path from $H$ to $F_j$ when the edge $z$ or a $z^*$-edge not derived from $z+v$ is used, then this  $z^*$-edge is always dropped

2.  on the path from $H$ to $F_i$ when the edge $v$ or a $v^*$-edge not derived from $z+v$ is used, then this  $v^*$-edge is always dropped

3. whenever on the path from $H$ to $F_j$ and $H$ to $F_i$ we have the same edges to do the reduction on, we go to the right or to the left on both paths

4. $F_j$ has a unique $z^*$-edge not derived from $z+v$ denoted by  $f(z)$ and $F_i$ has a unique $v^*$-edge not derived from $z+v$  denoted by $f(v)$, and we have $F_j-f(z)=F_i-f(v)$
\el

The proof of Lemma \ref{differ} can be seen by inspection and is left to the reader.

\medskip

\noindent{\it Proof of Lemma \ref{q_i}.} We prove by induction on $\alt(G)$, which the the maximum length of a path in $R_G^{\O}$ from $G$ to a leaf of $R_G^{\O}$,  that $P_i$ is the set of preceeding facets for $F_i$. Let $G_1, G_2, G_3$ be the left, right and middle children of $G$ and let $F_1, \ldots, F_k$, and $F_{k+1}, \ldots, F_l$ and $Q_1, \ldots, Q_z$ be their full dimensional leaves for which the statement holds. Thus, the set of preceeding facets of $F_i$, $i \in [k]$, is $P_i$ by inductive hypothesis. Assume that for some $k+1\leq i\leq l$ there is a preceeding facet $H$ of $F_i$ which is not in $P_i$. By the inductive hypothesis for $G_2$ and  the definition of preceeding facet, $H$ has to then be in $R_{G_3}^{\O}$. However, since $R_G^{\O}$ satisfies the weak embeddable property, we have that the only possible such facet is $b_G^{-1}(F_i)$,  when this is well defined. However, using that $R_G^{\O}$ satisfies the twosided embeddable property, we can then show that there exists a full dimensional leaf in $R_{G_1}^{\O}$ which contains $b_G^{-1}(F_i)$, thereby showing that $b_G^{-1}(F_i)$ is in $P_i$. Thus all preceeding facets of $F_i$ are in $P_i$. 

Next we need to show that all elements of $P_i$ are preceeding facets. First we observe that the elements of $P_i$, $k+1\leq i\leq l$,  must be leaves of $R_G^{\O}$. Indeed, if we are considering an element $F_j\cap F_i \in P_i$ with $k+1\leq j<i$, then it is a leaf of  $R_{G_2}^{\O}$ by inductive hypothesis. If  $F_j\cap F_i \in P_i$ for $j<k+1\leq i$, then we have $F_j$ in $R_{G_1}^{\O}$ and $F_i$ in $R_{G_2}^{\O}$, and the two graphs differ by exactly one edge. However, with Lemma \ref{differ} we can see that then $F_j \cap F_i$ is in $R_{G_3}^{\O}$, since we can obtain it by going towards $G_3$ from $G$ and then executing the operations as on the path to $F_j$ or $F_i$. 
 It is now clear that the elements of $P_i$  satisfy conditions 1. and 2. in Definition \ref{facet}. Assume that for some $k+1\leq i\leq l$ there is  $H\in P_i$ which is not a preceeding facet. Thus, because of the inductive hypothesis we have that $H \in R_{G_1}^{\O}$. However, if $H$ differs from $F_i$ by only missing an edge, then it can be seen that  the sequence of reductions used to obtain $H$ from $G$ vs the  sequence of reductions used to obtain $F_i$ from $G$ is different in that somewhere we need to go towards M for $H$ and towards $R$ to $F_i$. Thus, $H$ cannot belong to $R_{G_1}^{\O}$, completing the proof. 
\qed

\medskip

Before proceeding to the proof of Theorem \ref{thm-shell2} we introduce the strong embeddable property which, as its name suggests it is an extension of the weak embeddable property. We then see that $R_G^{\O}$ has this property and use it to prove 
Theorem \ref{thm-shell2}. The strong embeddable property will also be a basis for proofs of several nonnegativity results of reduced forms, including  the proof of a conjecture of Kirillov.

 \begin{definition} \label{strong} Let  $R_G$ posses the weak embeddable property. At a non-leaf $H$ of $R_G$ let $b_H$ be the bijection specified in 2 in Definition \ref{weak}. The reduction tree $R_G$ is said to have the (right) {\bf strong embeddable property} if the following statements are true:

1. if $b_H(Q_i)=F_{i_j}$,  so that $E(F_{i_j})=E(Q_i)\cup \{e\}$, then if in $R_{H_3}$ the preceeding facets of the full dimensional leaf $Q_i$ are  $Z_1, \ldots, Z_k$ (in the sense of Definition \ref{facet}), then $Z_1+e, \ldots, Z_k+e$  are preceeding facets of  $F_{i_j}$ in $R_{H_2}$  

2.  for $F_{i_j}$ as in 1, there are no  leaves in  $R_{H_2}$ which are preceeding  facets of  $F_{i_j}$  other than $Z_1+e, \ldots, Z_k+e$  
\end{definition}

Note, that if $R_H$ possesses the weak embeddable property then for a full dimensional leaf $L$ in $R_{H}$ which is also in $R_{H_2}$ there is exactly one preceeding facet of it belonging to $R_{H_3}$ if $L$ is in the image of $b_H$ and otherwise there is no preceeding facet of it belonging to $R_{H_3}$.

\begin{lemma} \label{st-embed} The reduction tree $R_G^{\O}$ has the strong embeddable property.
\end{lemma}

\noindent{\it Proof idea.}
 In light of Lemma \ref{q_i} strong embeddability of $R_G^{\O}$ is equivalent to \eqref{kz0} as explained below.
 At a non-leaf $H$ of $R_G^{\O}$ let $b_H$ be the bijection specified in 2 in Definition \ref{weak}.  Let $H_2$ and $H_3$ be the right and middle children of $H$ in $R_G^{\O}$. Let $F_{k+1}, \ldots, F_l$ be the full dimensional leaves of $R_{H_2}^{\O}$ and let 
 $Q_{1}, \ldots, Q_z$ be the full dimensional leaves of $R_{H_3}^{\O}$.   Let $b_H(Q_i)=F_{i_j}$,  so that $E(F_{i_j})=E(Q_i)\cup \{e_i\}$, $i \in [z]$.

For $k+2\leq i\leq l$,    let $$K_i=\{\{K_1^i, \ldots, K_{f(i)}^i\}\}=\{\{F_i \cap F_j \mid k+1\leq j<i, |E(F_i \cap F_j)|=|E(F_i)|-1\}\}.$$   

For $1\leq i\leq z$,    let $$Z_i=\{\{Z_1^i, \ldots, Z_{h(i)}^i\}\}=\{\{Q_i \cap Q_j \mid 1\leq j<i, |E(Q_i \cap Q_j)|=|E(Q_i)|-1\}\}.$$   

Then  we need to prove \be \label{kz0} \{\{K_1^{i_j}, \ldots, K_{f(i_j)}^{i_j}\}\}=\{\{Z_1^i+e_i, \ldots, Z_{h(i)}^i+e_i\}\}.\ee

Proving  \eqref{kz0}  can be accomplished by proving $K_{i_j} \subset Z_{i}+e$ and  $Z_{i}+e \subset K_{i_j}$ using case analysis and utilizing Lemma \ref{differ}.
\qed
 \medskip

\noindent \textit{Proof of Theorem \ref{thm-shell2}.}
We prove that \eqref{formal}  is the formal sum of the  set of the leaves of $R_G^{\O}$ by induction on $\alt(G)$. We know that $\sum_{i=1}^k \prod_{j=1}^{f(i)} (F_i+Q_j^i)$, $\sum_{i=k+1}^l \prod_{j=1}^{g(i)} (F_i+K_j^i)$, and $\sum_{i=1}^z \prod_{j=1}^{h(i)} (Q_i+Z_j^i)$
are the formal sums of the  set of the leaves of $R_{G_1}^{\O}$, $R_{G_2}^{\O}$ and $R_{G_3}^{\O}$, respectively, where the notation is as would be expected based on the statement of Theorem \ref{thm-shell2}. We  know by strong embeddability and Lemma \ref{q_i}  that if $b_G^{-1}(F_i)$, $k+1\leq i$,  is not well-defined then $\{Q_j^i\}_{i=1}^{f(i)}=\{K_j^i\}_{i=1}^{g(i)}$ and  if $b_G^{-1}(F_i)$ is well-defined then $\{Q_j^i\}_{i=1}^{f(i)}=\{K_j^i\}_{i=1}^{g(i)}\cup \{b_G^{-1}(F_i)\}$. But then  

\be \label{eq1} \sum_{i=1}^l \prod_{j=1}^{f(i)} (F_i+Q_j^i)=\sum_{i=1}^k \prod_{j=1}^{f(i)} (F_i+Q_j^i)+\sum_{i=k+1}^l \prod_{j=1}^{g(i)} (F_i+K_j^i)+\sum_{i=k+1}^l \prod_{j=1}^{g(i)} \chi(b_G^{-1}(F_i))(F_i+K_j^i)b_G^{-1}(F_i),\ee
where $\chi(b_G^{-1}(F_i))$ is $1$ is $b_G^{-1}(F_i)$ is well-defined, and $0$ otherwise. Note that

\be \label{bG} \sum_{i=k+1}^l \prod_{j=1}^{g(i)} \chi(b_G^{-1}(F_i))(F_i+K_j^i)b_G^{-1}(F_i)=\sum_{i=1}^z \prod_{j=1}^{g(b_G(i))} (b_G(Q_i)+K_j^{b_G(i)})Q_i,\ee where $b_G(i)$ is the index $k$ of $F_k$ to which $b_G(Q_i)$ is equal to.

The right hand side of \eqref{bG} is equal to:

\be \label{eq2} \sum_{i=1}^z \prod_{j=1}^{g(b_G(i))} ((b_G(Q_i) \cap Q_i)+(K_j^{b_G(i)}\cap Q_i))=\sum_{i=1}^z \prod_{j=1}^{h(i)} (Q_i+Z_j^i),\ee where the last equality holds by equation \eqref{kz0} stated in the proof of the strong embeddable property of $R_G^{\O}$ in Lemma \ref{st-embed}. 

Equations \eqref{eq1}, \eqref{bG} and \eqref{eq2} then imply

\be \sum_{i=1}^l \prod_{j=1}^{f(i)} (F_i+Q_j^i)= \sum_{i=1}^k \prod_{j=1}^{f(i)} (F_i+Q_j^i)+\sum_{i=k+1}^l \prod_{j=1}^{g(i)} (F_i+K_j^i)+\sum_{i=1}^z \prod_{j=1}^{h(i)} (Q_i+Z_j^i),\ee  completing the proof.

\qed

\section{Refining $h$-vectors of the  canonical triangulation of $\F_{\widetilde{G}}$ }
\label{sec:refining}


In this section we study a refinement of the $h$-polynomial of the canonical triangulation of flow polytopes.  

Consider the reduction tree $R_G^{\O}$ and let $F_i$ and $Q^i_j$ be as in  Theorem \ref{thm-shell2}.  By Theorem \ref{thm-shell2} each $Q^i_j$ appears in the reduction tree $R_G^{\O}$, and we assign a \textbf{weight} $w(Q_j^i)=\b_a$ to each $Q^i_j$, where the unique reduction on the path from $G$ to  $Q^i_j$ where we go to the middle child is performed on the edges $(a, c), (c, d)$.
By Theorem \ref{thm-shell2} all other not full dimensional simplices are obtained as intersections of a subset of the $Q^i_j$s, and we weight those intersections by the product of the weights of the $Q^i_j$ appearing in the intersection (note that this may or may not be the same as the product of $\b_i$'s associated to the sequence of reductions yielding the graph). 
Denote the weight of $G$ by $w(G)$. We set the weight of full dimensional leaves to be $1$. From what we just said, together with  Theorem \ref{thm-shell2}, it follows that: 

\bt

Let ${F_1}, \ldots,{F_l}$ be full dimensional leaves of $R_G^{\O}$ in depth-first search order. Let $$P_i:=\{\{Q_1^i, \ldots, Q_{f(i)}^i\}\}=\{\{F_i \cap F_j \mid 1\leq j<i, |E(F_i \cap F_j)|=|E(F_i)|-1\}\}.$$  Then

\be \label{huh2}  \sum_{i=1}^l \prod_{j=1}^{f(i)} (F_i+w(Q_j^i) Q_j^i),\ee
 \noindent is the formal sum of the  set of weighted  leaves of $R_G^{\O}$, where the product of graphs is their intersection,  and if $f(i)=0$ then we define $\prod_{j=1}^{f(i)} (F_i+w(Q_j^i)Q_j^i)=F_i$.
\et

 Let $\C$ be the abstract  simplicial complex obtained from $R_G^{\O}$, as in Theorem \ref{tri}. 
    Recall that $h(\C, \b)=\sum_{i=0}^d h_i \b^i,$ where using the shelling from Theorem \ref{thm-shell} we get that  $h_i$ is equal to the number of top dimensional simplices which attach on $i$ facets to the union of previous simplices in the shelling order. Equation \eqref{huh2} then suggests the following natural refinement of the $h$-vector of the  canonical triangulation of $\F_{\widetilde{G}}$ . 

\bd \label{frak}  For the   canonical triangulation of $\F_{\widetilde{G}}$  let the $h(\mathfrak{b})$-vector be the following refinement of the $h$-polynomial:
\be \label{def} h(\C, \mathfrak{b})=\sum_{i=1}^l\prod_{j=1}^{f(i)}w(Q_j^i).\ee Clearly, setting all $\b_i=\b$ we have $h(\C, \mathfrak{b})=h(\C, \b)$. Thus, $h(\C, \mathfrak{b})$ gives a refinement of the $h$-polynomial.
\ed

\bt \label{thm-o} Let $Q_G^{\O}(\mathfrak{b}; \x)$ be the reduced form in the subdivision algebra $\mt$ when we did the reductions in the specified order $\O$. Let $Q_G^{\O}(\mathfrak{b}-\bf{1})$ be the specialization of $Q_G^{\O}(\mathfrak{b}-{\bf 1}; \x)$ at $\x=(1, \ldots, 1)$.
Then  \be \label{o} Q_G^{\O}(\mathfrak{b}-{\bf1})=h(\C, \mathfrak{b}).\ee In particular, $Q_G^{\O}(\mathfrak{b}-\bf{1})$ has nonnegative integer coefficients.
\et

\proof 
We prove \eqref{o} by induction on $\alt(G)$.

If $\alt(G)=0$,  then $Q_G^{\O}(\mathfrak{b}-{\bf1})=1$ and $h(\C, \mathfrak{b})=1$, also.

Suppose \eqref{o} is true for all graphs $G$ with $\alt(G)<m$. Consider the graph $G$ with $\alt(G)=m>0$. Since there is a pair of alternating edges in $G$ we can perform a reduction on the edges $(i, j)$ and $(j, k)$, $i<j<k$, of  $G$  which come first in the order $\O$, obtaining the graphs $G_1$, $G_2$ and $G_3$, such that $\alt(G_1), \alt(G_2), \alt(G_3)<m$.   It follows then by definition that 
\be Q_G^{\O}(\mathfrak{b})=Q_{G_1}^{\O}(\mathfrak{b})+Q_{G_2}^{\O}(\mathfrak{b})+\b_i Q_{G_3}^{\O}(\mathfrak{b}).  \ee Since $\alt(G_1), \alt(G_2), \alt(G_3)<m$, it follows by inductive hypothesis that $Q_{G_i}^{\O}(\mathfrak{b}-1)=h(\C_i, \mathfrak{b})$, $i \in [3]$, where  $\C_i$ is the  canonical triangulation of $\F_{\widetilde{G}}$  of  ${\F_{\widetilde{G_i}}}$, $i \in [3]$. Next we show that 
\be \label{rec} h(\C, \mathfrak{b})=h(\C_1, \mathfrak{b})+h(\C_2, \mathfrak{b})+(\b_i-1) h(\C_3, \mathfrak{b}),  \ee which will conclude the proof of \eqref{o}.

Recall that by the definition of $h(\C, \mathfrak{b})$, we look at the shelling order $\F_{\wt{F_1}}, \ldots, \F_{\wt{F_l}}$ we obtained from reading off the full dimensional leaves of $R_G^{\O}$ in depth-first search order. By  Theorem \ref{thm-shell} we have that  $\F_{\wt{F_1}}, \ldots, \F_{\wt{F_f}}$ is  a shelling of $\F_{\wt{G_1}}$ read off from $R_{G_1}^{\O}$ and $\F_{\wt{F_{f+1}}}, \ldots, \F_{\wt{F_l}}$ is a shelling of $\F_{\wt{G_2}}$ read off from $R_{G_2}^{\O}$. Let $\F_{\wt{L_1}}, \ldots, \F_{\wt{L_s}}$ be a shelling of $\F_{G_3}$ read off from $R_{G_3}^{\O}$  in the same fashion. Note that 

 \be \label{c1} h(\C_1, \mathfrak{b})=\sum_{i=1}^f\prod_{j=1}^{f(i)}w(Q_j^i).\ee

Next we show that \be \label{c2} h(\C_2, \mathfrak{b})+(\b_i-1) h(\C_3, \mathfrak{b})=\sum_{i=f+1}^l\prod_{j=1}^{f(i)}w(Q_j^i).\ee
Equations \eqref{c1} and \eqref{c2} yield  \eqref{rec}.

Call $\prod_{j=1}^{f(i)}w(Q_j^i)$ the weight contribution of simplex $\F_{\wt{F_i}}$ to the $h( \mathfrak{b})$-polynomial. Note that a simplex $\F_{\wt{F_i}}$, $i \in \{f+1, \ldots, l\}$, contributes the same weight to   $h(\C_2, \mathfrak{b})$ and $h(\C, \mathfrak{b})$ if and only if  $\F_{\wt{F_i}} \cap(\F_{\wt{F_1}}\cup\cdots\cup \F_{\wt{F_f}})=\emptyset$. On the other hand using the strong embeddable property of $R_G^{\O}$,  the weight contribution of all simplices $\F_{\wt{F_i}}$, $i \in \{f+1, \ldots, l\}$ such that  $\F_{\wt{F_i}} \cap(\F_{\wt{F_1}}\cup\cdots\cup \F_{\wt{F_f}})\neq \emptyset$ is equal to $h(\C_3, \mathfrak{b})$ in $h(\C_2, \mathfrak{b})$ and $\b_i h(\C_3, \mathfrak{b})$ in $h(\C, \mathfrak{b})$, yielding \eqref{c2}.
\qed

Theorem \ref{thm-o} yields an alternative proof to \cite[Theorem 8]{h-poly1}.
Indeed,  by \cite[Lemma 5]{h-poly1} 
we have that  $Q_G^{\O}(\mathfrak{b}-\textbf{1})=Q_G(\b-1),$ when we set $\b_i=\b$. 
However, the initial proof of  \cite[Theorem 8]{h-poly1} 
is simpler then the above, building on much less knowledge.

An interesting special case of Theorem \ref{thm-o} to consider is when $G$ is the path graph $P_n=([n], \{(i, i+1) | i \in [n-1]\})$.  In this case the notion of weight $w(G)$ has an additional combinatorial interpretation. Before we proceed to state it, we note that the reduced form still depends on the order of reductions we use, and we keep to using the order $\O$ in the rest of this section.  Indeed, in the order $\O$ the leaf of $R_{P_5}^{\O}$  labeled by the graph $([5], \{(1, 5)\})$ is weighted by $\b_1^3$, whereas if we first reduce the edges $(1, 2)$ and $(2, 3)$ and then the edges $(3, 4)$ and $(4, 5)$, then it would be weighted by $\b_1^2 \b_3$.

Given a leaf $G$, let $((i_a, j_a), (j_a, k_a), M)$, $a\in [p]$, be all the reductions on the path from $P_n$ to $G$ in $R_{P_n}^{\O}$  where we go towards the middle child. Define $b(G)=\prod_{a=1}^p \b_{i_a}$ the \textbf{balance} of leaf $G$. The following theorem states that we can express $b(G)$ in terms of properties of $G$. 

\bt Given a leaf $G$ of $R_{P_n}^{\O}$, $$b(G)=\prod_{i=1}^{n-1}\b_i^{f_G(i)},$$ where $f_G(i)$ is equal to the number of (graph-)components of $G$ such that the shortest edge $e$ such that the component is entirely between the initial and end vertex of the edge $e$ has initial vertex $i$.
\et

\proof 


We prove by induction on $n$ that $b(G)=\prod_{i=1}^{n-1}\b_i^{f_G(i)}$. The base case is trivial. Assume it is true for all $P_m$, $m<n$. 

Consider $P_n$. Let $L_1, \ldots, L_k$ be the leaves of $R_{P_{n-1}}^{\O}$ in depth-first search order. By assumption, $b(L_j)=\prod_{i=1}^{n-2}\b_i^{f_{L_j(i)}}.$ Consider a leaf $L$ of $R_{P_{n}}^{\O}$. Since the leaves of $R_{P_{n}}^{\O}$ are the leaves of $R_{L_1+(n-1,n)}^{\O},\cdots,R_{L_k+(n-1, n)}^{\O}$, we can assume that $L$ is a leaf of $R_{L_z+(n-1,n)}^{\O}$, $z \in [k]$. Then we have that \be \label{mult}b_{P_n}(L)=b_{P_{n-1}}(L_z)\times b_{L_z+(n-1,n)}(L),\ee where we indexed $b$  to clarify within which reduction tree we are. Combining \eqref{mult} with the inductive hypothesis yields our desired result.
\qed

\bc Given a leaf $G$ of $R_{P_n}^{\O}$ with $n-2$ edges, $$w(G)=\prod_{i=1}^{n-1}\b_i^{f_G(i)},$$ where $f_G(i)$ is equal to the number of (graph-)components of $G$ such that the shortest edge $e$ such that the component is entirely between the initial and end vertex of the edge $e$ has initial vertex $i$.

\ec

\proof This is immediate, since for a leaf $G$ with $n-2$ edges the weight $w(G)$ is defined to be equal to the balance $b(G)$.
\qed

It appears to be true in general that if $G$ is a leaf in $R_{P_n}^{\O}$, then $w(G)=b(G)$.  We leave this investigation to the interested reader.

\section{The weak and strong embeddable properties of partial reduction trees}
\label{sec:embed}

In this section we define partial reduction trees and show how to extend the previous theorems to them. 
Reduction trees in the algebras $\yb$ and $\myb$ can be considered as partial reduction trees  in the sense of this section, so the results presented below can be used for studying reduced forms in $\yb$ and $\myb$.

\bd \label{partial} A {\bf partial reduction tree} of  the reduction tree $R_G$ is a connected rooted  (at $G$) subtree of $R_G$ such that if a vertex has a left or middle or right child, then it has all three. We denote a partial reduction tree by $R_G^{p}$. 
\ed

\begin{definition} \label{weak-p} A partial reduction tree $R_G^p$ is said to have the (right) {\bf weak embeddable property} if one of the following is true for every vertex $H$ of $R_G^p$:  

1. $H$ is a leaf of $R_G^p$

2. the middle child of $H$ is $H_3$ and the right child of $H$ is $H_2$, satisfying that there is a map $b_H$ from the full dimensional leaves of the subtree $R_{H_3}^p$ of $R_G^p$ (leaves with $|E(H_3)|$ number of edges) into  the full dimensional leaves of the subtree $R_{H_2}^p$ (leaves with $|E(H_2)|$ number of edges) of $R_G^p$ such that if $b_H(L)=L'$, then $E(L')=E(L) \cup \{e\}$ with $e \not \in E(L)$. Moreover, if $L'$ is a full dimensional leaf of $R_{H_2}^p$ such that there is a leaf $L$ of  $R_{H_3}^p$ with the property that $E(L')=E(L) \cup \{e\}$ with $e \not \in E(L)$, then $L$ is in the image of $b_H$ and $b_H$ is a bijection from the full dimensional leaves of  $R_{H_3}^p$  onto its image.  Moreover, there is a unique $L$ in $R_{H_3}$ such that  $E(L')=E(L) \cup \{e\}$ with $e \not \in E(L)$.
\end{definition}

\bd \label{q-facet} Given a partial reduction tree  $R_G^p$ and a full dimensional leaf $L$ of it, we say that a leaf $H$ of  $R_G^p$  is a {\bf preceeding  facet} of $L$ if 

1. $H$ is before $L$ in the depth first search order of  the leaves of   $R_G^p$ 

2. $E(H)\subset E(L)$ and $|E(H)|=|E(L)|-1$

3. the unique path in  $R_G^p$ from $L$ to $H$ consists of several up steps followed by several down steps, so that the first of the down steps is a Middle reduction. 

\ed

\begin{definition} \label{strong-p} Let  $R_G^p$ posses the weak embeddable property. At a non-leaf $H$ of $R_G^p$ let $b_H$ be the bijection specified in 2 in Definition \ref{weak-p}. The partial reduction tree $R_G^p$ is said to have the (right) {\bf strong embeddable property} if the following statements are true:

1. if $b_H(Q_i)=F_{i_j}$,  so that $E(F_{i_j})=E(Q_i)\cup \{e\}$, then if in $R_{H_3}^p$ the preceeding facets of the full dimensional leaf $Q_i$ are  $Z_1, \ldots, Z_k$ (in the sense of Definition \ref{q-facet}), then $Z_1+e, \ldots, Z_k+e$  are preceeding facets of  $F_{i_j}$ in $R_{H_2}^p$  

2.  for $F_{i_j}$ as in 1, there are no  leaves in  $R_{H_2}^p$ which are preceeding  facets of  $F_{i_j}$  other than $Z_1+e, \ldots, Z_k+e$  
\end{definition}

Note that if $R_G^p$ possesses the weak embeddable property then for a full dimensional leaf $L$ in $R_{H}^p$ which is also in $R_{H_2}^p$ there is exactly one preceeding facet of it belonging to $R_{H_3}^p$ if $L$ is in the image of $b_H$ and otherwise there is no preceeding facet of it belonging to $R_{H_3}^p$.


\bd Given a partial reduction tree  $R_G^p$ with the strong embeddable property define the  $h(\mathfrak{b})$-polynomial for it as follows:

\be \label{b} h(R_G^p, \mathfrak{b})=\sum_{L} p(L),\ee where the sum runs over all  full dimensional leaves $L$ of $R_G^p$ and $p(L)=\prod_{F} w(F),$ where the product  is over the  preceeding facets $F$ of  $L$ and $w(F)$ is the weight of $F$ as defined in Section \ref{sec:refining}. The empty product is defined to be equal to $1$.

If we specialize by setting $\beta_i=\beta$ for all $i \in [n]$, then we get
\be h(R_G^p, \beta)=\sum_{i=0}^{\infty}s_i \beta^i,\ee where $s_i$ is the number of full dimensional leaves $L$ of $R_G^p$ such that there are exactly $i$ preceeding facets of it. \ed 

Note that if we take $R_G^p$ to be the reduction tree $R_G^{\O}$ then the  $h(\mathfrak{b})$-polynomial in \eqref{b}  agrees with the $h(\mathfrak{b})$-polynomial from Definition \ref{frak} and specializes to the usual $h$-polynomial.

The following result is a culmination of the insight of the above definitions. It generalizes Theorem \ref{thm-o} and \cite[Theorem 8]{h-poly1}.

\bt  \label{q-h} Given a partial reduction tree  $R_G^p$ with the strong embeddable property we have that \be \label{p-red} Q_G^p(\mathfrak{b}-1)=h(R_G^p, \mathfrak{b}),\ee
where $Q_G^p(\mathfrak{b}; \x)=\sum_L \x(L) \mathfrak{b}(L)$, where the sum is over all  leaves of $R_G^p$,  $\x(L)=\prod_{(i,j) \in L} x_{ij}$ and $ \mathfrak{b}(L)=\prod_{j=1}^z \beta_{i_j}$, where on the path from $G$ to $L$ we went towards the middle $z$ times, and the reductions where we went towards the middle had the minimal vertex of the first edge be $i_1, \ldots, i_z$.
We denote $Q_G^p(\mathfrak{b}-1)=Q_G^p(\mathfrak{b}-1; \bf{1})$.
\et


\proof We prove Theorem \ref{q-h} by induction on $\alt(G)$. Since  $R_G^p$ has the strong embeddable property, so do  $R_{G_1}^p$, $R_{G_2}^p$  and $R_{G_3}^p$, where $G_1, G_2, G_3$ are as in \eqref{graphs1} after we performed reduction on the edges $(i, j)$ and $(j, k)$ of $G$. By definition,  
\be Q_G^p(\mathfrak{b})=Q_{G_1}^p(\mathfrak{b})+Q_{G_2}^p(\mathfrak{b})+\beta_i Q_{G_3}^p(\mathfrak{b}), \ee thus to prove \eqref{p-red} we need to prove that

\be \label{h} h(R_G^p, \mathfrak{b})=h(R_{G_1}^p, \mathfrak{b})+h(R_{G_2}^p, \mathfrak{b})+(\beta_i-1) h(R_{G_3}^p,\mathfrak{b}), \ee holds, since by induction  $Q_{G_i}^p(\mathfrak{b}-1)=h(R_{G_i}^p, \mathfrak{b})$, for $i \in [3]$. Equation \eqref{h} follows by definition, since the strong embeddable property ensures  that the contribution of $R_{G_1}^p$ to $h(R_G^p, \mathfrak{b})$ is exactly $h(R_{G_1}^p, \mathfrak{b})$ and  the contribution of $R_{G_2}^p$ to $h(R_G^p, \mathfrak{b})$ is exactly $h(R_{G_2}^p, \mathfrak{b})+(\beta_i-1) h(R_{G_3}^p,\mathfrak{b})$, as explained in the following. The full dimensional leaves of $R_{G_2}^p$ which are not in the image of the map $b_G$  contribute  $h(R_{G_2}^p, \mathfrak{b})-h(R_{G_3}^p,\mathfrak{b})$ to  $h(R_G^p, \mathfrak{b})$  and  full dimensional leaves of $R_{G_2}^p$ which are in the image of  $b_G$  contribute $\beta_i h(R_{G_3}^p,\mathfrak{b})$ to  $h(R_G^p, \mathfrak{b})$,  where we multiply by $\beta_i$ since if a full dimensional leaves of $R_{G_2}^p$ is in the image of   $b_G$, then other than the preceeding facets of it in $R_{G_2}^p$, it has one additional preceeding facet, namely its preimage under $b_G$.
\qed

\section{Solving   Kirillov's Conjecture 7} 
\label{sec:conjecture7}

 In this section we use the techniques developed in Section \ref{sec:embed} to prove  \cite[Conjecture 7]{k2} for a special family of reduction trees, namely those which posses the extra strong embeddable property. We also demonstrate  via counterexamples that   \cite[Conjecture 7]{k2}  fails in general.   We recall the conjecture here for convenience.

\bd Given a graph $G$ on the vertex set $[n]$, denote by $Q_G^{\mt}(\mathfrak{b}, t)$ the specialization of a particular reduced form $Q_G^{\mt}(\mathfrak{b}, \x)$ when  $x_{ij} = 1$, if $(i,j) \neq (1, n)$, and $x_{1,n} = t$. The polynomial $Q_G^{\mt}(\mathfrak{b}, t)$  depends on the order of reductions performed. If we set $\beta_i=\beta$ for all $i \in [n-2]$, then we denote $Q_G^{\mt}(\mathfrak{b}, t)$  by $Q_G^{\t}(\beta, t)$. 
\ed

\begin{conjecture} \cite[Conjecture 7 (A)]{k2} Let $n\geq 4$ and write \be \label{7a} Q_{K_n}^{\t}(\beta, t)=\sum_{k=0}^{2n-6} (1+\beta)^k c_{k,n}(t).\ee Then $c_{k,n}(t) \in \mathbb{Z}_{\geq 0}[t].$
\end{conjecture}

Note that since the reduced forms $Q_G^{\mt}(\mathfrak{b}, t)$  and $Q_G^{\t}(\beta, t)$ depend on the order of reductions performed, not just on the initial monomial determined by $G$, Kirillov's conjectures says that for any particular reduced form there is an expansion of the given form.  

\bd \label{t} Given a partial reduction tree  $R_G^p$ with the strong embeddable property define the generalized $h(\mathfrak{b}, t)$-polynomial for it as follows:

\be \label{b} h(R_G^p, \mathfrak{b}, t)=\sum_{L} p(L, t),\ee where the sum runs over all  full dimensional leaves $L$ of $R_G^p$ and $p(L,t)=w_t(L)\prod_{F} w(F),$ where the product  is over the  preceeding facets $F$ of  $L$ and $w_t(L)=t^l$ if $L$ has exactly $l$ edges $(1, n)$. The empty product is defined to be equal to $1$.
 
\ed

\bd A partial reduction tree $R_G^p$ is said to have the {\bf extra strong embeddable property} if it has the strong embeddable property and in addition for every non-leaf vertex $H$ the map $b_H$ maps a graphs with a given number of edges $(1, n)$ to graphs with the same number of edges $(1, n)$.
\ed


The following result is a culmination of the insight of the above definitions.

\bt  \label{q-ht} Given a partial reduction tree  $R_G^p$ with the extra strong embeddable property we have that \be \label{q-p-red} Q_G^p(\mathfrak{b}-1,t)=h(R_G^p, \mathfrak{b}, t).\ee

\et


\proof We prove Theorem \ref{q-ht} by induction on $\alt(G)$. Since  $R_G^p$ has the strong embeddable property, so do  $R_{G_1}^p$, $R_{G_2}^p$  and $R_{G_3}^p$, where $G_1, G_2, G_3$ are as in \eqref{graphs1} after we performed reduction on the edges $(i, j)$ and $(j, k)$ of $G$. By definition,  
\be Q_G^p(\mathfrak{b},t)=Q_{G_1}^p(\mathfrak{b},t)+Q_{G_2}^p(\mathfrak{b},t)+\beta_i Q_{G_3}^p(\mathfrak{b},t), \ee thus to prove \eqref{q-p-red} we need to prove that

\be \label{eq-q-h} h(R_G^p, \mathfrak{b},t)=h(R_{G_1}^p, \mathfrak{b},t)+h(R_{G_2}^p, \mathfrak{b},t)+(\beta_i-1) h(R_{G_3}^p,\mathfrak{b},t), \ee holds, since by induction  $Q_{G_i}^p(\mathfrak{b}-1,t)=h(R_{G_i}^p, \mathfrak{b},t)$, for $i \in [3]$. Equation \eqref{q-h} follows by definition, since the extra strong embeddable property ensures  that the contribution of $R_{G_1}^p$ to $h(R_G^p, \mathfrak{b},t)$ is exactly $h(R_{G_1}^p, \mathfrak{b},t)$ and  the contribution of $R_{G_2}^p$ to $h(R_G^p, \mathfrak{b},t)$ is exactly $h(R_{G_2}^p, \mathfrak{b},t)+(\beta_i-1) h(R_{G_3}^p,\mathfrak{b},t)$, as explained in the following. The full dimensional leaves of $R_{G_2}^p$ which are not in the image of the map $b_G$  contribute  $h(R_{G_2}^p, \mathfrak{b},t)-h(R_{G_3}^p,\mathfrak{b},t)$ to  $h(R_G^p, \mathfrak{b},t)$  and  full dimensional leaves of $R_{G_2}^p$ which are in the image of  $b_G$  contribute $\beta_i h(R_{G_3}^p,\mathfrak{b},t)$ to  $h(R_G^p, \mathfrak{b},t)$,  where we multiply by $\beta_i$ since if a full dimensional leaves of $R_{G_2}^p$ is in the image of   $b_G$, then other than the preceeding facets of it in $R_{G_2}^p$, it has one additional preceeding facet, namely its preimage under $b_G$.
\qed

\bt \label{7} Suppose that the reduced form  $Q_G^{\mt}(\mathfrak{b}, t)$ was obtained through a reduction tree with the extra strong embeddable property and write 
 \be \label{p7} Q_{G}^{\mt}(\mathfrak{b}-1, t)=\sum_{k=0}^{\infty} \sum_{I : |I|=k} p(I) c_{I}(t),\ee where the sum is over all multisets $I$ with elements in $[n-2]$ and with cardinality $k$, and $p_k(I)=\prod_{i\in I}\beta_i$. Then  $c_{I}(t) \in \mathbb{Z}_{\geq 0}[t].$
Thus, 
 \cite[Conjecture 7 (A)]{k2} holds for any graph (not just complete) if the corresponding reduction tree has the extra  strong embeddable property. 
\et

\proof By Theorem \ref{q-ht} we have that $Q_{G}^{\mt}(\mathfrak{b}-1, t)=h(R_G, \mathfrak{b}, t),$ where $R_G$ is a reduction tree with the extra strong embeddable property. Together with Definition \ref{t} this implies the statement of Theorem \ref{7} immediately.
\qed

Since Theorem \ref{7} proves   \cite[Conjecture 7 (A)]{k2} only if the reduction tree has the extra strong embeddable property, it raises the question of what happens otherwise. In Figure \ref{fig:ce} we present the smallest counterexample to  \cite[Conjecture 7 (A)]{k2}. As can be seen the reduction tree in Figure \ref{fig:ce} does not have the extra strong embeddable property. This counterexample was constructed by a program kindly written by  Leonid Chindelevitch. The same program found many other counterexamples  to  \cite[Conjecture 7 (A)]{k2}.

\begin{figure}
\begin{center}
\includegraphics[scale=.9]{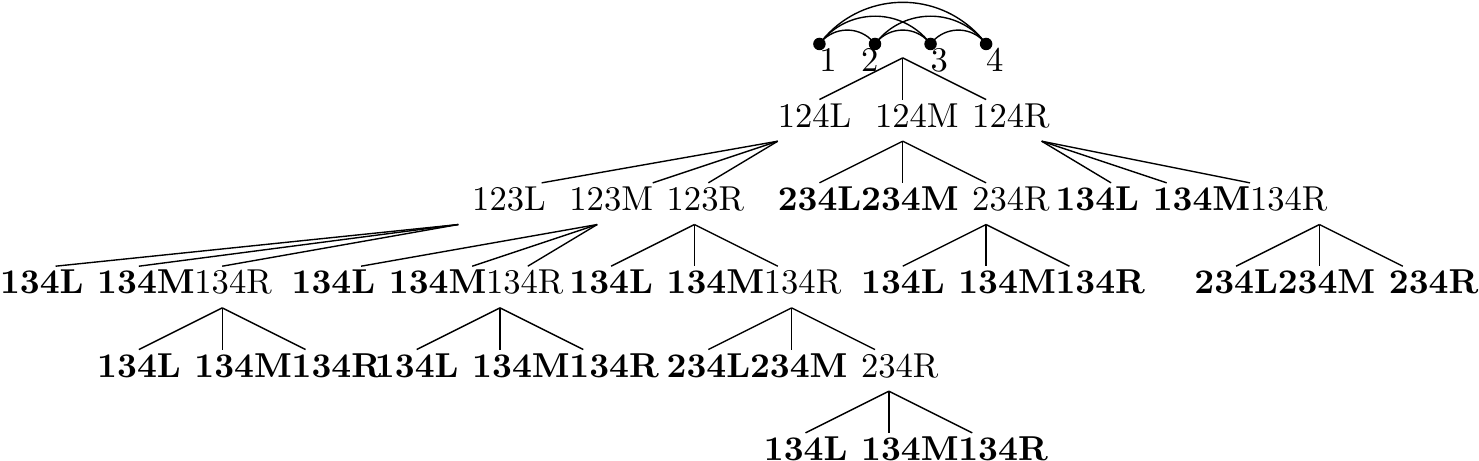}
 \caption{A reduction tree of the complete graph $K_4$. The label $ijkX$ at a node, where $1\leq i<j<k\leq 4$ and $X\in \{L, M, R\}$, specifies that the graph with this label was obtained by performing reduction $((i,j),(j,k), X)$ on its parent. When $ijkX$ is in boldface, it indicates that the corresponding graph is alternating. For this particular reduction tree we have $Q_{G}^{\t}(\b-1, t)= \b^0 t^4+ \b^1(-t^2+4t^3+2t^4)+\b^2(t^2+2t^3+t^4)$, contradicting  \cite[Conjecture 7 (A)]{k2}.}
 \label{fig:ce}
 \end{center}
\end{figure}

  \section*{Acknowledgements}
The author is grateful to Leonid Chindelevitch for his interest, time and effort invested in writing the program which found the counterexample appearing in Section \ref{sec:conjecture7}. 
The author thanks Federico Ardila, Louis Billera and Ed Swartz for several inspiring discussions regarding this work.  The author also thanks Alejandro H. Morales for many  interesting conversations about flow polytopes over the course of the years.

\newpage
\section{Appendix}
\label{sec:app}

The purpose of this section is to provide a proof of Theorem \ref{tri} as well as several auxiliary results  used in the text. While we found it easier to define flow polytopes $\F_{\widetilde{G}}$ in Definition \ref{leaf} as a convex hull of vertices, in this section we will prefer their definition as intersection of halfspaces. As such, we now proceed to introduce some new notation and then  prove Theorem \ref{tri}.

\bd \label{2} Fix a graph $G$ on the vertex set $[n+1]$. Fix an order  on the edges of the complete graph $K_{n+1}$: $e_1, \ldots, e_{{n+1 \choose 2}}$.  Let $x_1, \ldots, x_k$ be the set of {\bf base variables}, by which we mean a set of variables in which everything will be expressed. Let $\c=(c_1, c_2, \ldots)$ be an infinite  vector with finitely many nonzero entries whose coordinates $c_i$, $i \in [m]$, satisfy the following two properties: 

\begin{itemize}

\item $c_i$, $i \in [m]$,  is a linear combination of the base variables,

\item $c_i=0$, if $i=a{n+1 \choose 2}+b$, $0<b\leq {n+1 \choose 2}$, and the multiplicity of edge $e_b$ is less than $a+1$ in $G$. 

\end{itemize}
\ed

In other words, we think of the $i$th coordinate $c_i$  of $\c$ as corresponding to a possible edge in $G$, namely to  the  $(a+1)$st edge $e_b$, if $i=a{n+1 \choose 2}+b$, $0<b\leq {n+1 \choose 2}$. Then, the above requirements say that if an edge is not in $G$, then set the corresponding variable to $0$, and otherwise to a linear combination of the base variables. Given an edge $e=(i,j)$ in $G$, we also write $c(e)$ for the corresponding variable in $\c$. Namely, if $(i, j)=e_b$ in the ordering of the edges of $K_{n+1}$ and $e$ is the $(a+1)st$ edge $(i,j)$ in $G$, then $c(e)=c_{a{n+1 \choose 2}+b}.$

\bd \label{fp} The {\bf flow polytope} $\F_G(\c)$ (with base variables $x_1, \ldots, x_k$) is defined as   $x_{i}\geq 0$, $i \in k$, $c_i\geq 0$, $i \in \Z_{\geq 0}$, and

  $$1=\sum_{e \in E(G), \i(e)=1}c(e)= \sum_{e \in E(G),  \fin(e)=n+1}c(e),$$
  
\noindent  and for $2\leq i\leq n$
  
  $$\sum_{e \in E(G), \fin(e)=i}c(e)= \sum_{e \in E(G), \i(e)=i}c(e).$$
  
\ed

Note that if we order the multiset of edges $E(G)=\{\{e_1, \ldots, e_l\}\}$ and take variables $x_{e_i}$, $i \in [l]$, to be the base variables in Definition \ref{2}, and let $c(e_i)=x_{e_i}$, $i \in [l]$, then $\F_G(\c)$ of Definition \ref{2} is the usual way to define flow polytopes $\F_G$.

 

Recall that given a graph $G$ on the vertex set $[n+1]$  containing edges $(i, j)$ and $(j, k)$, $i<j<k$, performing the {\bf reduction} on these edges of $G$ yields three graphs on the vertex set $[n+1]$:
  \begin{eqnarray} \label{graphs10}
E(G_1)&=&E(G)\backslash \{(j, k)\} \cup \{(i, k)\}, \nonumber \\
E(G_2)&=&E(G)\backslash \{(i, j)\} \cup \{(i, k)\},  \nonumber \\
E(G_3)&=&E(G)\backslash \{(i, j), (j, k)\} \cup \{(i, k)\}.
\end{eqnarray}
  
 Suppose that the edge $(i,j)$ of $G$ involved in the reduction is the $d$th among the edges $(i,j)\in E(G)$, the edge $(j,k)$ of $G$ involved in the reduction is the $f$th among the edges $(j,k)\in E(G)$ and   there are $a\geq 0$ edges $(i, k)$ present in the graph $G$. Consider the flow polytope $\F_{\tilde{G}}[\c]$ as in Definition \ref{2}. Define $\c^{G_1}, \c^{G_2}, \c^{G_3}$ to agree with $\c$ on all coordinates except on the coordinates corresponding to the $d$th edge $(i, j)$, $f$th edge $(j, k)$ and $(a+1)$st edge $(i,k)$. Call these edges $g_1, g_2, g_3$ for simplicity. Then, when $c(g_1)\geq c(g_2)$ we can write $c^{G_1}(g_1)=c(g_1)-c(g_2)$, $c^{G_1}(g_2)=0$, $c^{G_1}(g_3)=c(g_2)$. When $c(g_2)\geq c(g_1)$ we can write $c^{G_2}(g_1)=0$, $c^{G_2}(g_2)=c(g_2)-c(g_1)$, $c^{G_2}(g_3)=c(g_1)$. Finally, when $c(g_1)= c(g_2)$ we can write $c^{G_3}(g_1)=0$, $c^{G_3}(g_2)=0$, $c^{G_1}(g_3)=c(g_1)=c(g_2)$. At times if clarity requires we denote $\c^{G_1}, \c^{G_2}, \c^{G_3}$ by $\c^{G_1, G}, \c^{G_2, G}, \c^{G_3, G}$ to emphasize that the reduction tree is rooted at $G$ and the base variables correspond to the edges of $G$.

Using the above, we can restate Lemma \ref{red} as follows; only our notation has changed.

\bl{ \cite[Proposition 1]{m-prod},\cite[Proposition 4.1]{mm}, \cite{p, S}}
 \label{red1} Given a graph $G$ on the vertex set $[n+1]$ and   $(i, j), (j, k) \in E(G)$,   for some $i<j<k$,  we have  $$\F_{{\tilde{G}}}(\c)=\F_{{\tilde{G}_1}}(\c^{G_1}) \bigcup \F_{{\tilde{G}_2}}(\c^{G_2}),  \F_{{\tilde{G}_1}}(\c^{G_1}) \bigcap \F_{{\tilde{G}_2}}(\c^{G_2})=\F_{{\tilde{G}_3}}(\c^{G_3}) \text{ and } \F_{{\tilde{G}_1}}(\c^{G_1})^\circ \bigcap \F_{{\tilde{G}_2}}^\circ(\c^{G_2})=\emptyset,$$

 \noindent where $\F_{{\tilde{G}}}(\c)$, $\F_{{\tilde{G}_1}}(\c^{G_1})$, $\F_{{\tilde{G}_2}}(\c^{G_2})$ are of the same dimension $d-1$, $\F_{{\tilde{G}_3}}(\c^{G_3})$ is $d-2$ dimensional,  and  $\mathcal{P}^\circ$ denotes the interior of $\mathcal{P}$.
\el

\bd \label{4} Let $H$ be a graph labeling a node of the reduction tree $R_G$. Let the unique path from $G$ to $H$ be (in terms of the graphs on the nodes) $G-I_1-\cdots-I_p-H$. Definition \ref{3} constructs $\c^{I_1, G}$. Successively applying the rules given in Definition \ref{3}, while keeping the base variables those corresponding to the edges of $G$, we obtain the vector $\c^{H, G}$.

\ed

\noindent \textbf{Intersection of flow polytopes as intersection of  graphs.} We show that if we use the special order $\O$ to reduce the graphs we consider, then in a sense (made precise below) we can think of intersections of two flow polytopes as intersection of  graphs. Such a property is in general unexpected, and highlights the special choice of our reduction order $\O$.

\bd \label{g_1capg_2} Given vectors $\c^{G_1, G}$ and $\c^{G_2, G}$ as in Definition \ref{4}, we define vector $\c_{G_1\cap G_2, G}$  as follows. Considering 
$\c^{G_1, G}$ and $\c^{G_2, G}$ as vectors expressed in the variables $x_1, \ldots, x_k$ satisfying constraints as in  Definition \ref{4}, let  $\c_{G_1\cap G_2, G}$  be the vector we obtain if we require  additionally that  $\c^{G_1, G}=\c^{G_2, G}$, putting additional constraints on the variables $x_1, \ldots, x_k$.  That is, treating  each coordinate $(\c^{G_1, G})_i$,$(\c^{G_2, G})_i$,  $ i\geq 1$, as an expression in $x_1, \ldots, x_k$, let $C$ be the set of conditions on $x_1, \ldots, x_k$ from Definition \ref{4} arising because of the vectors $\c^{G_1, G}$ and $\c^{G_2, G}$  together with the conditions $(\c^{G_1, G})_i=(\c^{G_2, G})_i$,  $ i\geq 1$. Then 
$(\c_{G_1\cap G_2, G})_i=(\c^{G_1, G})_i|_{C}$, where $(\c^{G_1, G})_i|_{C}$ is equal to $(\c^{G_1, G})_i$ with the conditions $C$ satisfied.
\ed

The purpose of Definition \ref{g_1capg_2} is to express the intersection of two flow polytopes corresponding to leaves $G_1$ and $G_2$ of $R_G^{\O}$, as stated in the following theorem.

\bt \label{coord} Let $G_1$ and $G_2$ be two leaves of $R_G^{\O}$. Then 

\be \label{int}  \F_{\widetilde{G_1}}({\c^{G_1, G}})\cap \F_{\widetilde{G_2}}(\c^{G_2, G})=\F_{\widetilde{G_1\cap G_2}}(\c_{G_1\cap G_2, G}).\ee  Moreover, if $G_1\cap G_2$ is  in $R_G^{\O}$, then $\c^{G_1\cap G_2, G}=\c_{G_1\cap G_2, G}.$ (Note that  $G_1\cap G_2$ need not be in $R_G^{\O}$.)
\et

Before proving Theorem \ref{coord} we need to provide several auxiliary results. 


\bt \label{c} If $G_1$ and $G_2$ are leaves of $R_G^{\O}$, then the vector $\c_{G_1\cap G_2, G}$  can be obtained from $\c^{G_1, G}$ by setting $c^{G_1, G}(e)=0$ for all $e \not \in G_1 \cap G_2$. In other words, the conditions on the base variables posed by $\c^{G_1, G}$ and $c^{G_1, G}(e)=0$ for all $e \not \in G_1 \cap G_2$ are equivalent to the conditions on the base variables posed by $\c^{G_2, G}$ and $c^{G_2, G}(e)=0$ for all $e \not \in G_1 \cap G_2$ and both of these are equivalent to the conditions on the base variables posed by  $\c^{G_1, G}$ and $\c^{G_2, G}$ together with the conditions $(\c^{G_1, G})_i=(\c^{G_2, G})_i$,  $ i\geq 1$. 
\et

Before proceeding to the proof of Theorem \ref{c}, we prove the following special case of it: 

\bp  \label{p-c}Theorem \ref{c} holds for graphs $G$ with the property that  there exists a vertex $v \in V(G)$ such that all edges of $G$ are incident to $v$ and $v$ has $k>0$ incoming edges and one outgoing edge. 
\ep 

\proof We prove Proposition \ref{p-c} by induction on $k$. When $k=1$, it is trivial to check the statement. Assume the statement is true for all $m<k$. Let $G$ be a graph with $k$ incoming edges into $v$ and one outgoing edge. Let $e$ be the outgoing edge and $e'$ be the lowest incoming edge into $v$. Consider $G'=G-e'$. By the inductive hypothesis the statement of Proposition \ref{p-c} holds for it.  Let $F_1, \ldots, F_{k}$ be the full dimensional leafs of   $R_{G'}^{\O}$ in depth-first search order.   Note that since $e'$ was the lowest edge and we are doing reduction in order $\O$, it follows that the reduction tree of $G$ can be obtained from the reduction tree of $G'$ by adding the edge $e'$ to the leaves of $R_{G'}^{\O}$ and reducing where necessary. It is easy to see that there is only one leaf of $R_{G'}^{\O}$, namely $F_{k}$, where adding the edge $e'$ makes the leaf a nonalternating graph. Let the leaves of $R_{F_{k}}^{\O}$ be $F_{k}'$, $Q'$  and $F_{k+1}'$ in depth-first search order. Now consider two leaves $G_1$ and $G_2$ of $R_{G}^{\O}$. If  $G_1, G_2 \not \in \{F_{k}', F_{k+1}',Q'\}$ then $G_1-e'$ and $G_2-e'$  are leaves of  $R_{G'}^{\O}$ and we can conclude the statement of Theorem \ref{c} for them by inductive hypothesis. If $G_1, G_2 \in \{F_{k}', F_{k+1}',Q'\}$, then it is easy to check   the statement of Theorem \ref{c} for them directly (since we can effectively consider $F_k$ as the root of the reduction tree). Finally, if  $G_1 \not \in \{F_{k}', F_{k+1}',Q'\}$ and $G_2  \in \{F_{k}', F_{k+1}',Q'\}$ we consider the three cases depending on whether $G_2=F_{k}',  G_2=F_{k+1}'$ or $G_2=Q'$. In all three cases it suffices to remember that the edge we obtain from performing the reduction on $e'$ and $e$ is not part of any  leaf of $R_G^{\O}$ other than $F_{k}', F_{k+1}',Q'$, and that $e'$ is an edge of all the other leaves of  $R_G^{\O}$. Using these two facts and that $F_{k}', F_{k+1}',Q'$ are the children of $F_k$,  the statements can be derived readily from the fact that the statement of Proposition \ref{p-c} holds for $G_1-e'$ and $F_k$. 

\qed

\noindent \textit{Proof  of Theorem \ref{c}.}   We prove Theorem \ref{c} by induction on $\alt(G)$. 

\medskip 
\noindent {\bf Base of induction.} $\alt(G)=0$. In  this case $\F_{\widetilde{G}}$ is a simplex and the statement is trivial.

\medskip 
\noindent {\bf Inductive hypothesis.} Statement true for all $G$,   $\alt(G)<m$, $m>0$.

\medskip 
\noindent {\bf Inductive step.} Consider $G$ with $\alt(G)=m$.
Look at the last vertex $v$ of $G$ which is nonalternating. Let $e$ be the lowest edge outgoing from $v$. Then $\alt(G-e)<m$, so the statement is true for it. From here we will prove that it is also true for $G$. 

Let $G_1$ and $G_2$ be two leaves of $R_{G-e}^{\O}$. By the inductive hypothesis Theorem \ref{c} holds for them. 
Consider two leaves $H_1$ and $H_2$ of $R_{G}^{\O}$.  There are two possible cases:

{\it Case 1.} $H_1$ and $H_2$ are leaves of $R_{G_1+e}^{\O}$, where $G_1$ is a leaf of $R_{G-e}^{\O}$.

{\it Case 2.} $H_1$ and $H_2$ are leaves of $R_{G_1+e}^{\O}$ and $R_{G_2+e}^{\O}$, respectively,  where $G_1$ and $G_2$ are distinct  leaves of $R_{G-e}^{\O}$.

In Case 1, we can use coordinates $\c^{H_1, G_1+e}$ and $\c^{H_2, G_1+e}$, since whatever functions of the original edges of $G-e$  the base variables corresponding to the edges of $G_1$ are, the combinations remain untouched as we proceed with reductions in  $R_{G_1+e}^{\O}$. Since $G_1$ is a leaf of  $R_{G-e}^{\O}$, either $G_1+e$ is alternating, in which case we are done, or it has exactly one nonalternating vertex $v$ with one outgoing edge $e$ and some incoming edges $e_1, \ldots, e_k$, $k \geq 1$. In this case we can use Proposition \ref{p-c} directly to prove that Theorem \ref{c} holds for $H_1$ and $H_2$. 

In Case 2 we 
 need to consider cases based on whether $G_1+e$ and $G_2+e$ are both alternating, both nonalternating, or one is alternating and one is nonalternating. In all these subcases, the statement of Theorem  \ref{c} for $H_1$ and $H_2$ follows from 
 Proposition \ref{p-c} together with the fact  that we distinguish edges $(i, j)$ based on how they were obtained as explained in Definition \ref{edges} -- that is, if an edge $(i, j)$ was obtained by doing a reduction on two edges $e_1$ and $e_2$ and another edge $(i, j)$ was obtained by  doing a reduction on two edges $e_3$ and $e_4$ with $\{e_1, e_2\}\neq \{e_3, e_4\}$, then these two edges $(i, j)$ are considered different, and the vectors $\c$ take this into account. \qed



Now we are ready to prove  Theorem \ref{coord}.

 \medskip

\noindent {\it Proof of Theorem \ref{coord}.}   By the definition of  $\c_{G_1\cap G_2, G}$  it follows that  
\be \F_{\widetilde{G_1\cap G_2}}(\c_{G_1\cap G_2, G}) \subset \F_{\widetilde{G_1}}({\c^{G_1, G}})\cap \F_{\widetilde{G_2}}(\c^{G_2, G}).\ee

To show \be \label{sub} \F_{\widetilde{G_1}}({\c^{G_1, G}})\cap \F_{\widetilde{G_2}}(\c^{G_2, G}) \subset \F_{\widetilde{G_1\cap G_2}}(\c_{G_1\cap G_2, G})\ee
consider a point $\p \in \F_{\widetilde{G_1}}({\c^{G_1, G}})\cap \F_{\widetilde{G_2}}(\c^{G_2, G}).$ Since $\p \in \F_{\widetilde{G_1}}({\c^{G_1, G}})$, it follows that all nonzero coordinates of $\p$ lie on the edges of $G_1$ and since $\p \in \F_{\widetilde{G_2}}(\c^{G_2, G})$, all nonzero coordinates of $\p$ lie on the edges of $G_2$. Thinking of the coordinates of $\p$ as corresponding to edges and expressing it with respect to the base variables corresponding to the edges of $\widetilde{G}$, we can conclude that  $\p$ is a particular evaluation of the variable coordinate 
vector $\c^{G_1, G}$ with the constraint that $\c^{G_1, G}(e)=0$ for all $e \not \in G_1$. Similarly, $\p$ is a particular evaluation of the variable coordinate 
vector $\c^{G_2, G}$ with the constraint that $\c^{G_2, G}(e)=0$ for all $e \not \in G_2$. Thus, using the meaning given to $\c_{G_1\cap G_2, G}$ in Theorem \ref{c}, equation \eqref{sub} follows. 

Next we note that if for two leaves $G_1$ and $G_2$ of  $R_G^{\O}$ the intersection $G_1\cap G_2$ is also in $R_G^{\O}$, then the sequence of reductions used to obtain  $G_1\cap G_2$ can be obtained from the sequence of reductions used to obtain $G_1$ by going towards the middle in some of the reductions and accordingly deleting some followng reductions. It follows that  $\c^{G_1\cap G_2, G}$  can be obtained from $\c^{G_1, G}$ by setting $c^{G_1, G}(e)=0$ for all $e \not \in G_1 \cap G_2$. Therefore, by Theorem \ref{c} we have that $\c^{G_1\cap G_2, G}=\c_{G_1\cap G_2, G}$.
\qed

\medskip
Now we are ready to prove Theorem \ref{tri}. For convenience we repeat its statement here.
\medskip

\noindent \textbf{Theorem 2.} \textit{The simplices corresponding to the  full dimensional leaves of  $R_{G}^{\O}$ induce a triangulation; that is, the intersection of any two of them is a face of both. Moreover, the simplices corresponding to all leaves of  $R_{G}^{\O}$ are part of this triangulation.}
\proof By Theorem \ref{coord} the intersection  $\F_{\widetilde{G_1}}({\c^{G_1, G}})\cap \F_{\widetilde{G_2}}(\c^{G_2, G})$ is $\F_{\widetilde{G_1\cap G_2}}(\c_{G_1\cap G_2, G})$. It follows readily from Theorem \ref{c}  that $\F_{\widetilde{G_1\cap G_2}}(\c_{G_1\cap G_2, G})$ is a face of both $\F_{\widetilde{G_1}}({\c^{G_1, G}})$ and  $\F_{\widetilde{G_2}}(\c^{G_2, G})$. \qed

The following lemma is important for determining the dimension of $ \F_{\widetilde{G_1\cap G_2}}(\c_{G_1\cap G_2, G})$ using Theorem \ref{coord}.

\bl As  polynomials in the base variables corresponding to the edges of G, $\c_{G_1\cap G_2, G}(e)$ for $e \in G_1\cap G_2$ are linearly independent.
\el

\proof By Theorem \ref{c}, $\c_{G_1\cap G_2, G}=\c^{G_1, G}$ subject to the constraints $c^{G_1, G}(e)=0$ for all $e \not \in G_1 \cap G_2$. We claim that the coordinate polynomials in $\c^{G_1, G}$ corresponding to the edges of $G_1$ are linearly independent. This would   imply that $\c_{G_1\cap G_2, G}(e)$ for $e \in G_1\cap G_2$ are also linearly independent.

To see that the coordinate polynomials in $\c^{G_1, G}$ corresponding to the edges of $G_1$ are linearly independent observe that this is true of $\c^{G, G}$, which basically says that the base variables corresponding to the edges of $G$ are distinct.  Looking at the path in $\R_G^{\O}$ from $G$ to $G_1$ we can prove the claim by induction on $\alt$. Observe that in each step we take a pair of coordinate polynomials $(p_1, p_2)$ and replace them by either $(p_1, p_2-p_1)$, $(p_2, p_1-p_2)$ or by a single polynomial $p_1$. Clearly then the resulting coordinate polynomials are still linearly independent.
\qed

\bc \label{dimension} The dimension of  $\F_{\widetilde{G_1\cap G_2}}(\c_{G_1\cap G_2, G})$ is $|E(G_1\cap G_2)|+|V(G_1\cap G_2)|-1$.
\ec

\proof The dimension of a flow polytope of $G$ is $|E(G)|-|V(G)|+1$ \cite{BV2}. Since $|E(\widetilde{G_1\cap G_2})|=|E(G_1\cap G_2)|+2|V(G_1\cap G_2)|$, $ |V(\widetilde{G_1\cap G_2})|=|V(G_1\cap G_2)|+2$, the result follows.

\qed
\bt  \label{subsets} Let $G_1, G_2$ and $G_3$ be three leaves of $R_G^{\O}$ so that $G_1\cap G_2 \subset G_1 \cap G_3$. Then 

\be \label{int1} \F_{\widetilde{G_1\cap G_2}}(\c_{G_1\cap G_2, G}) \subset  \F_{\widetilde{G_1\cap G_3}}(\c_{G_1\cap G_3, G}).\ee   
\et

\proof Recall that  by Theorem \ref{c},  $\c_{G_1\cap G_2, G}$  can be obtained from $\c^{G_1, G}$ by setting $c^{G_1, G}(e)=0$ for all $e \not \in G_1 \cap G_2$ and  $\c_{G_1\cap G_3, G}$  can be obtained from $\c^{G_1, G}$ by setting $c^{G_1, G}(e)=0$ for all $e \not \in G_1 \cap G_3$. Thus,  $\c_{G_1\cap G_2, G}$  can be obtained from $\c_{G_1\cap G_3, G}$ by setting $\c_{G_1\cap G_3, G}(e)=0$ for all $e  \in (G_1\cap G_3)-(G_1 \cap G_2)$. 
Thus given $\p \in \F_{\widetilde{G_1\cap G_2}}(\c_{G_1\cap G_2, G})$ it follows that $\p \in \F_{\widetilde{G_1\cap G_3}}(\c_{G_1\cap G_3, G})$ proving \eqref{int1}.
\qed

\bibliography{biblio-kir}

\begin{thebibliography}{10}

\bibitem{BV2}
W.~Baldoni and M.~Vergne.
\newblock {K}ostant partitions functions and flow polytopes.
\newblock {\em Transform. Groups}, 13(3-4):447--469, 2008.

\bibitem{kosh}
V.~I. Danilov, Karzanov A.V., and G.~A. Koshevoy.
\newblock Coherent fans in the space of flows in framed graphs.
\newblock {\em DMTSC proc., FPSAC 2012 Nagoya, Japan}, pages 483--494, 2012.

\bibitem{fk}
S.~Fomin and A.~N. Kirillov.
\newblock Quadratic algebras, dunkl elements, and schubert calculus.
\newblock {\em Advances in Geometry}, (172).

\bibitem{Kir}
A.~N. Kirillov.
\newblock Ubiquity of {K}ostka polynomials.
\newblock {\em arXiv}, math.QA, Dec 1999.

\bibitem{k2}
A.N. Kirillov.
\newblock On some combinatorial and algebraic properties of {D}unkl elements.
\newblock {\em RIMS preprint}, 2012.

\bibitem{k2014}
A.N. Kirillov.
\newblock On some quadratic algebras, {D}unkl elements, schubert, grothendieck,
  tutte and reduced polynomials.
\newblock {\em RIMS preprint}, 2014.

\bibitem{m-prod}
K.~M\'esz\'aros.
\newblock Product formulas for volumes of flow polytopes.
\newblock {\em Proc. Amer. Math. Soc., to appear}.

\bibitem{m2}
K.~M\'esz\'aros.
\newblock Root polytopes, triangulations, and the subdivision algebra, {I}.
\newblock {\em Trans. Amer. Math. Soc.}, 363(8):4359--4382, 2011.

\bibitem{m1}
K.~M\'esz\'aros.
\newblock Root polytopes, triangulations, and the subdivision algebra, {I}{I}.
\newblock {\em Trans. Amer. Math. Soc.}, 363(11):6111--6141, 2011.

\bibitem{h-poly1}
K.~M\'esz\'aros.
\newblock $h$-polynomials via reduced forms.
\newblock 2014.
\newblock \href{http://arxiv.org/abs/math/1018992}{arXiv:1407.2685}.

\bibitem{mm}
K.~M\'esz\'aros and A.~H. Morales.
\newblock Flow polytopes of signed graphs and the {K}ostant partition function.
\newblock {\em International Mathematical Research Notices, to appear}.

\bibitem{p}
A.~Postnikov, 2010.
\newblock personal communication.

\bibitem{S}
R.P. Stanley.
\newblock Acyclic flow polytopes and {K}ostant's partition function.
\newblock Conference transparencies,
  \url{http://math.mit.edu/~rstan/trans.html}, 2000.

\end{thebibliography}
\bibliographystyle{plain}

\end{document}